\documentclass[11pt]{article} 

\usepackage{amsmath,amssymb,amsfonts,float,theorem,enumitem} 
\usepackage[all]{xy}

\def\zr#1{\bibitem{#1} A.~Zinger,}

\numberwithin{equation}{section} 
\newtheorem{thm}{Theorem}[section]
\newtheorem{prp}[thm]{Proposition} 
\newtheorem{lmm}[thm]{Lemma}  
\newtheorem{dfn}[thm]{Definition}
\newtheorem{crl}[thm]{Corollary}
\theorembodyfont{\rmfamily}

\newtheorem{eg}{Example}

\def\BE#1{\begin{equation}\label{#1}}  
\def\EE{\end{equation}}

\def\ti#1{\widetilde{#1}} 
\def\e_ref#1{(\ref{#1})} 
\def\lan{\langle} \def\ran{\rangle} 
\def\lr#1{\lan#1\ran} 
\def\blr#1{\big\lan#1\big\ran}
\def\bllrr#1{\big\lan\!\big\lan{#1}\big\ran\!\big\ran}
\def\lra{\longrightarrow}
\def\sf#1{\textsf{#1}}
\def\ov#1{\overline{#1}}

\def\Lra{\Longrightarrow} \def\Llra{\Longleftrightarrow}

\def\al{\alpha}
\def\be{\beta}
\def\ga{\gamma}

\def\ep{\epsilon} 
\def\io{\iota}
\def\ka{\kappa}
\def\na{\nabla}
\def\om{\omega}
\def\si{\sigma}

\def\ups{\upsilon}

\def\vph{\varphi}
\def\vr{\varrho}
\def\ze{\zeta}

\def\A{\mathcal A}
\def\C{\mathbb C}

\def\bC{\mathbf C}
\def\cD{\mathcal D}
\def\cE{\mathcal E}

\def\H{\mathcal H}

\def\cJ{\mathcal J}

\def\K{\mathcal K}

\def\bL{\mathbb L}
\def\M{\mathfrak M}

\def\O{\mathcal O}
\def\Q{\mathbb Q}

\def\P{\mathbb P}
\def\R{\mathbb R}

\def\T{\mathcal T}
\def\U{\mathcal U}
\def\fU{\mathfrak U}
\def\X{\mathfrak X}
\def\W{\mathcal W}

\def\Z{\mathbb Z}

\def\De{\Delta}
\def\Ga{\Gamma}
\def\La{\Lambda}

\def\Si{\Sigma}

\def\a{\mathbf a}

\def\f{\mathbf f}
\def\fI{\mathfrak i}
\def\fJ{\mathfrak j}
\def\u{\mathbf u}

\def\i{\infty}
\def\part{\partial}
\def\dbar{\bar\partial}
\def\eset{\emptyset} 

\def\Ann{\textnormal{Ann}}
\def\Bev{\textnormal{\bf ev}}

\def\cok{\textnormal{cok}}

\def\cok{\textnormal{cok}}

\def\eff{\textnormal{eff}}

\def\ev{\textnormal{ev}}

\def\GL{\textnormal{GL}}
\def\h{\textnormal{h}}

\def\Hom{\textnormal{Hom}}
\def\id{\textnormal{id}}
\def\Id{\textnormal{Id}}
\def\Im{\textnormal{Im}}
\def\ind{\textnormal{ind}}

\def\PD{\textnormal{PD}}

\def\reg{\textnormal{reg}}
\def\Re{\textnormal{Re}}
\def\Res{\textnormal{Res}}
\def\rk{\textnormal{rk}}
\def\sgn{\textnormal{sgn}}
\def\sing{\textnormal{sing}}
\def\std{\textnormal{std}}
\def\st{\textnormal{s.t.}}
\def\v{\textnormal{v}}

\begin{document}

\title{A Comparison Theorem for Gromov-Witten Invariants\\
in the Symplectic Category}

\author{Aleksey Zinger\thanks{Partially supported by a Sloan fellowship 
and DMS Grant 0604874}}
\date{\today}

\maketitle

\begin{abstract}
\noindent
We exploit the geometric approach to the virtual fundamental class,
due to Fukaya-Ono and Li-Tian, 
to compare Gromov-Witten invariants of 
a symplectic manifold and a symplectic submanifold whenever 
all constrained stable maps to the former are contained in the latter to first order.
Various special cases of the comparison theorem in this paper 
have long been used in the algebraic category;
some of them have also appeared in the symplectic setting.
Combined with the inherent flexibility of the symplectic category,
the main theorem leads to a confirmation 
of Pandharipande's Gopakumar-Vafa prediction for GW-invariants
of Fano classes in 6-dimensional symplectic manifolds.
The proof of the main theorem uses deformations of the Cauchy-Riemann
equation that respect the submanifold and 
Carleman Similarity Principle for solutions of perturbed Cauchy-Riemann
equations.
In a forthcoming paper, we apply a similar approach to relative Gromov-Witten invariants
and the absolute/relative correspondence in genus~$0$.
\end{abstract}

\tableofcontents

\section{Introduction}
\label{intro_sec}

\noindent
Gromov-Witten invariants are certain counts of pseudo-holomorphic curves
in symplectic manifolds that play prominent roles in symplectic topology,
algebraic geometry, and string theory.
These are usually rational numbers, and their precise relations with
some sort of integer enumerative counts of curves are rarely clear.
However, it is well-known that  genus~$0$ GW-invariants of
Fano manifolds are precisely counts of rational curves;
this observation is key to enumerating rational curves in projective space 
in \cite[Section~5]{KoM} and \cite[Section~10]{RT}.
String theory predicts an amazing integral structure for GW-invariants of Calabi-Yau
threefolds. These predictions originate in \cite{AsMo}, \cite{GV1}, and~\cite{GV2}
and are extended to all threefolds in~\cite{P1}.\\

\noindent
GW-invariants of a symplectic manifold~$X$ are obtained by evaluating 
natural coholomogy classes on
the virtual fundamental class (\sf{VFC}) of the space of stable 
$J$-holomorphic maps to~$X$.
The main statement of this paper, Theorem~\ref{main_thm}, compares 
GW-invariants counting stable maps meeting specified constraints 
in the ambient manifold 
with analogous counts of such maps to a submanifold containing 
the images of all such constrained maps to first order.
It leads immediately to Corollary~\ref{LeP_crl}, 
which in a way is a succinct re-formulation of the main conclusion of~\cite{LeP},
and with a bit more work to Theorem~\ref{FanoGV_thm}, which confirms 
the ``Fano case" of the Gopakumar-Vafa prediction of \cite[Section~0.2]{P1}.
Theorem~\ref{main_thm} is obtained by deforming the Cauchy-Riemann equation
in two stages so that the first stage respects the submanifold.
Carleman Similarity Principle is used to take advantage of properties
of solutions of Cauchy-Riemann equations that are preserved by
a large class of perturbations of the equations.
In a forthcoming paper \cite{divisorGWs}, we 
will apply similar geometric principles to study relative 
GW-invariants and the absolute/relative correspondence in genus~$0$ with 
applications to birational geometry in the spirit of Hu-Li-Ruan
(\cite{HuLtR}, \cite{HuR}, \cite{LtR}) and McDuff (\cite{Mc}).\\

\noindent
The author would like to thank R.~Pandharipande for bringing the 
``Fano case" of the Gopakumar-Vafa prediction to the author's attention,
D.~McDuff for detailed comments and suggestions on an earlier version of this paper, 
and T.~Graber, T.-J.~Li, D.~Maulik, and Y.~Ruan for related discussions.

\subsection{A comparison theorem for GW-invariants}
\label{AbsGW_subs}

\noindent
We will denote by $\bar\Z^+$ the set of non-negative integers.
Let $(X,\om)$ be a compact symplectic manifold.
If $g\!\in\!\bar\Z^+$, $S$ is a finite set, $\be\!\in\!H_2(X;\Z)$, and 
$J$ is an $\om$-tame\footnote{an almost complex structure on $(X,\om)$ is
\sf{$\om$-tame} if $\om(v,Jv)>0$ for all $v\!\in\!TX$ with $v\!\neq\!0$} 
almost complex structure on $X$, 
denote by $\ov\M_{g,S}(X,\be;J)$ the moduli space of equivalence classes of
stable $S$-marked genus~$g$ degree~$\be$ $J$-holomorphic maps to~$X$.
For each $j\!\in\!S$, there is a well-defined evaluation map
\BE{evmap_e} \ev_j\!: \ov\M_{g,S}(X,\be;J)\lra X.\EE
As standard in GW-theory, we will denote by
$$\psi_j\in H^2\big(\ov\M_{g,S}(X,\be;J)\big)$$
the first chern class of the universal cotangent line bundle for the $j$-th marked point.
The space $\ov\M_{g,S}(X,\be;J)$ carries a natural VFC,
which is independent of $J$ and will be denoted by  $[\ov\M_{g,S}(X,\be)]^{vir}$.
If the (real) dimension of~$X$ is $2n$, then
\BE{virdim_e}\dim\big[\ov\M_{g,S}(X,\be)\big]^{vir}
=\dim_{g,S}(X,\be) \equiv
2\big(\lr{c_1(TX),\be}+(n\!-\!3)(1\!-\!g)+|S|\big).\EE
If $J$ is regular\footnote{an almost complex structure $J$ is \sf{genus~$0$ regular}
if for every $J$-holomorphic map $u\!:\Si\!\lra\!X$, 
where $\Si$ is a tree of Riemann spheres, 
the linearization $D_{J;u}$ of the $\dbar_J$-operator at~$u$ is surjective},
then  $\ov\M_{0,S}(X,\be;J)$ is a topological manifold with a preferred choice 
of orientation and 
$$\big[\ov\M_{0,S}(X,\be)\big]^{vir}=\big[\ov\M_{0,S}(X,\be;J)\big].$$\\

\noindent
If $a_j\!\in\!\bar\Z^+$ and $\ka_j\!\in\!H_*(X;\Z)$ for each $j\!\in\!S$, let
\BE{GWdfn_e}\big((\tau_{a_j}\ka_j)_{j\in S}\big)_{g,\be}^X \equiv 
\bigg\lan\prod_{j\in S}
\big(\psi_j^{a_j}\ev_j^*(\PD_X\ka_j)\big),\big[\ov\M_{g,S}(X,\be)\big]^{vir}\bigg\ran,\EE
where $\PD_X\ka_j\!\in\!H^*(X;\Z)$ is the Poincare dual of $\ka_j$ 
in~$X$.\footnote{In the descriptions of Sections~\ref{semipos_subs} 
and~\ref{mainpf_subs}, 
$\big[\ov\M_{g,S}(X,\be)\big]^{vir}$ is a homology class in 
an arbitrarily small neighborhood of $\ov\M_{g,S}(X,\be;J)$ 
in the space of equivalence classes of $L^p_1$-maps to~$X$; 
there are well-defined evaluation maps~$\ev_j$ and cohomology classes~$\psi_j$ 
on this space as well.}
In order to avoid any sign ambiguities, we  define the number in~\e_ref{GWdfn_e} 
to be~$0$ if the dimension of~$\ka_j$ is odd for some~$j$.
By~\e_ref{virdim_e}, this number is zero unless 
\BE{dimcond_e}\sum_{j\in S}(2a_j+2n-\dim\ka_j)=\dim_{g,S}(X,\be).\EE
The number \e_ref{GWdfn_e} can be expressed as an integral on a ``smaller" moduli
space as follows.
Choose cobordism representatives $f_j\!:M_j\!\lra\!X$ for $\ka_j$, 
with $j\!\in\!S$.\footnote{We can assume that this is possible, since each $\ka_j$
can be replaced by a multiple for our purposes.}
Let
\BE{Mfdfn_e}\ov\M_{g,\f}(X,\be;J)=\big\{\big([u],(w_j)_{j\in S}\big)
\in \ov\M_{g,S}\big(X,\be;J\big)\!\times\!\prod_{j\in S}\!M_j\!:
\ev_j([u])\!=\!f_j(w_j)~\forall\,j\!\in\!S\big\}.\EE
The space $\ov\M_{g,\f}(X,\be;J)$ of constrained stable maps
also carries a virtual fundamental class and
$$\big((\tau_{a_j}\ka_j)_{j\in S}\big)_{g,\be}^X =
\bigg\lan\prod_{j\in S}\psi_j^{a_j},\big[\ov\M_{g,\f}(X,\be;J)\big]^{vir}\bigg\ran.$$
The subject of this section is a reduction of this GW-invariant of $X$ 
to a combination of GW-invariants for its submanifolds.

\begin{dfn}\label{propinter_dfn}
Let $Y$ be a submanifold of $X$.
A smooth map $f\!:M\!\lra\!X$ intersects~$Y$ \sf{properly} if 
$f^{-1}(Y)\!\subset\!M$ is a smooth orientable even-dimensional submanifold of~$M$ and 
$$d_wf\big(T_w\big(f^{-1}(Y)\big)\big)= d_w(TM)\cap T_{f(w)}Y$$
for every $w\!\in\!f^{-1}(Y)$.
\end{dfn}

\noindent
If $f\!:M\!\lra\!X$ intersects $Y\!\subset\!X$ transversally
and $M$, $X$, and $Y$ are orientable
of even total dimension, then $f$ intersects~$Y$ properly.
However, a proper intersection need not be transverse.
For example, any two real lines in~$\R^n$ intersect properly, 
but not transversally if $n\!\ge\!3$.
Two curves that are tangent to each other do not intersect properly.\\

\noindent
If $f\!:M\!\lra\!X$ intersects $Y\!\subset\!X$ properly and $NY\!\lra\!Y$ is
the normal bundle of~$Y$ in~$X$, the homomorphisms
$$d_w^{NY}f\!: T_wM\lra N_{f(w)}Y, \quad v\lra d_wf(v)+T_{f(w)}Y,
\qquad w\!\in\!f^{-1}(Y),$$
have constant rank; the kernel of $d_w^{NY}f$ is $T_w(f^{-1}(Y))$.
If $M$, $X$, and $Y$ are oriented, an orientation on $f^{-1}(Y)$ then induces 
an orientation on the vector bundle
$$N^fY\equiv f^*NY\big/(\Im\, d^{NY}f)\lra f^{-1}(Y).$$
Note that
\BE{rkNf_e} \rk\,N^fY=\big(\dim X-\dim M\big)-\big(\dim Y-\dim f^{-1}(Y)\big).\EE\\

\noindent
Let $Y$ be a compact symplectic submanifold of $X$ and 
$$\io_{Y*}\!: H_*(Y;\Z)\lra H_*(X;\Z)$$
the homomorphism induced by the inclusion $\io_Y\!:Y\!\lra\!X$.
If $\be_Y\!\in\!H_2(Y;\Z)$ and $J$ is an $\om$-tame almost complex structure on~$X$ 
which preserves $TY\!\subset\!TX|_Y$, then~$\io_Y$ induces an embedding
$$ \ov\M_{g,S}(Y,\be_Y;J)\hookrightarrow \ov\M_{g,S}(X,\io_{Y*}\be_Y;J).$$
If $f_j\!:M_j\!\lra\!X$, $j\!\in\!S$, are smooth maps as above, let
$$\ov\M_{g,\f}(Y,\be_Y;J)=\big\{\big([u],(w_j)_{j\in S}\big)
\in \ov\M_{g,\f}(X,\io_{Y*}\be_Y;J)\!:\, [u]\!\in\!\ov\M_{g,S}(Y,\be_Y;J)\big\}.$$
If in addition $u\!:\Si_u\!\lra\!Y$ is a $J$-holomorphic map from a nodal Riemann
surface (see Section~\ref{NRS_subs}), 
let $\H_u$ denote the space of deformations of the complex structure on~$\Si_u$.
The linearization of the $\dbar_J$-operator for maps to~$X$,
$$D_{J;u}^X\!:\H_u\oplus
L^p_1(\Si_u;u^*TX)\lra L^p(\Si_u;T^*\Si_u^{0,1}\!\otimes\!_{\C}u^*TX),
\quad p\!>\!2,$$
induces a generalized Cauchy-Riemann operator 
$$D_{J;u}^{NY}\!:  L^p_1(\Si_u;u^*NY)\lra 
L^p(\Si_u;T^*\Si_u^{0,1}\!\otimes\!_{\C}u^*NY).$$
For each $j\!\in\!S$, define
$$\ti\ev_j\!:\ker D_{J;u}^{NY}\lra N_{z_j(u)}Y
\qquad\hbox{by}\quad \xi\lra \xi(z_j(u))+T_{z_j(u)}Y,$$
where $z_j(u)\!\in\!\Si_u$ is the $j$-th marked point;
this homomorphism is the composition of the differential of the 
evaluation map~\e_ref{evmap_e} with the projection to the normal bundle.

\begin{thm}
\label{main_thm}
Suppose $(X,\om)$ is a compact symplectic $2n$-manifold, $g\!\in\!\bar\Z^+$,
$S$ is a finite set, $\be\!\in\!H_2(X;\Z)$, $a_j\!\in\!\bar\Z^+$ for each $j\!\in\!S$,
and $f_j\!:M_j\!\lra\!X$ is a cobordism representative for $\ka_j\!\in\!H_*(X;\Z)$
for each $j\!\in\!S$.
If $J$ is an $\om$-tame almost complex structure on~$X$,
$Y$ is a compact almost complex submanifold of~$(X,J)$, and
$\be_Y\!\in\!H_2(Y;\Z)$ are such that 
\begin{enumerate}[label=(\alph*)]
\item $\io_{Y*}(\be_Y)\!=\!\be$ and $f_j$ intersects $Y$ properly for each 
$j\!\in\!S$;
\item for every $([u],(w_j)_{j\in S})\!\in\!\ov\M_{g,\f}(Y,\be_Y;J)$, the homomorphism
\BE{mainthm_e}
\ker(D_{J;u}^{NY})\lra \bigoplus_{j\in S}N_{f_j(w_j)}^{f_j}Y, \qquad
\xi\lra \big(\ti\ev_j(\xi)+(\Im\,d_{w_j}f_j)\big)_{j\in S}\,,\EE
is an isomorphism,
\end{enumerate}  
then
\begin{enumerate}[label=(\arabic*)]
\item the space $\ov\M_{g,\f}(Y,\be_Y;J)$ carries a natural VFC 
(dependent on the orientations of $f_j^{-1}(Y)$) with
\BE{Ydim_e}\begin{split}
\dim \big[\ov\M_{g,\f}(Y,\be_Y;J)\big]^{vir}
&=\dim_{g,S}(X,\be)-\sum_{j\in S}\big(2n\!-\!\dim\ka_j\big)\\
&\qquad
+\sum_{j\in S}\rk\,N^{f_j}Y-2\big(\lr{c_1(NY),\be_Y}+\rk_{\C}NY\cdot(1\!-\!g)\big);
\end{split}\EE
\item the vector spaces $\cok(D_{J;u}^{NY})$ form a natural oriented vector orbi-bundle
$$\cok\big(D_J^{NY}\big) \lra \ov\M_{g,\f}(Y,\be_Y;J)$$
with
\BE{cokrk_e}
\rk_{\R} \cok(D_{J;u}^{NY})=
\sum_{j\in S}\rk\,N^{f_j}Y-2\big(\lr{c_1(NY),\be_Y}+\rk_{\C}NY\cdot(1\!-\!g)\big);\EE
\item $\ov\M_{g,\f}(Y,\be_Y;J)$ is a union of connected components of $\ov\M_{g,\f}(X,\be;J)$
and its contribution to the number~\e_ref{GWdfn_e} is given~by
\BE{mainthm_e2}\bC_{g,\f}(Y,\be_Y)=
\bigg\lan e\big(\cok(D_J^{NY})\big)\prod_{j\in S}\psi_j^{a_j},
\big[\ov\M_{g,\f}(Y,\be_Y;J)\big]^{vir}\bigg\ran.\EE\\
\end{enumerate} 
\end{thm}

\begin{eg}\label{CY_eg}
Suppose $(X,J)$ is a Calabi-Yau $3$-fold and $Y\!\subset\!X$ 
is a smooth isolated rational curve with $NY\!\approx\!\O(-1)\!\oplus\!\O(-1)$.
We can then apply Theorem~\ref{main_thm} with $S\!=\!\eset$, $g\!=\!0$, and 
$\be\!=\!d\io_{Y*}([Y])$ for any $d\!\in\!\Z^+$.
The assumption on the normal bundle implies that $\ker(D_{J;u}^{NY})$
is trivial and thus Condition~(b) is satisfied.
The right-hand side of~\e_ref{mainthm_e2} is then the famous multiple-cover 
contribution of~$1/d^3$ (\cite{AsMo}, \cite[Section~27.5]{MirSym}, \cite{Voisin}).
\end{eg}

\begin{eg}\label{inY_eg}
If the image of each map $f_j$ in Theorem~1.2 lies in $Y$,
the second part of Condition~(a) is automatically satisfied.
Condition~(b) is equivalent to the homomorphisms
$$\bigoplus_{j\in S}\ti\ev_j\!:
\ker(D_{J;u}^{NY})\lra \bigoplus_{j\in S}N_{z_j(u)}Y,
\qquad ([u],w)\in \ov\M_{g,\f}(X,\be;J),$$
being isomorphisms.
For example, this is the case if $X\!=\!\P^n$, $Y\!=\!\P^1\!\subset\!X$, $S\!=\!\{1,2\}$,
$g\!=\!0$, $\be\!=\!\io_{Y*}([Y])$ is the homology class of a line, $a_1,a_2\!=\!0$, 
and $f_1,f_2\!:pt\!\lra\!Y$ are maps to two distinct points.
In this particular case, 
$$\ov\M_{0,\f}(X,\be;J)=\ov\M_{0,\f}(Y,\be_Y;J),$$
where $\be_Y\!=\![Y]$, and $\cok(D_J^{NY})$ is the zero vector bundle.
Thus, 
$$\big(pt,pt\big)_{0,\be}^{\P^n}=
\big((\tau_{a_j}\ka_j)_{j\in S}\big)_{0,\be}^X
=\bC_{0,\f}(Y,\be_Y)
=\,^{\pm}\big|\ov\M_{g,\f}(Y,\be_Y;J)\big|
=\big(pt,pt\big)_{0,\be_Y}^{\P^1}
=1,$$
as expected.\footnote{\label{inY_ft}This is the number of lines through $2$ points in $\P^n$.
In this particular case, each operator $D_{J;u}^{NY}$ is $\C$-linear and 
its zero-dimensional kernel is positively oriented.
In general, this need not be the case; see \cite[Sections~9,10]{LeP}
for explicit sign computations.}
\end{eg}

\begin{eg}\label{transversetoY_eg}
If each map $f_j$ in Theorem~1.2 is transverse to $Y$,
the second part of Condition~(a) is again automatically satisfied.
Condition~(b) is equivalent to the injectivity of the operators $D_{J;u}^{NY}$
whenever $([u],w)\!\in\!\ov\M_{g,\f}(X,\be;J)$.
For example, this is the case if $X$ is the blowup of~$\P^n$, with $n\!\ge\!2$,
at a point, $Y\!\approx\!\P^{n-1}$ is the exceptional divisor,
$S\!=\!\{1,2\}$, $g\!=\!0$, $\be_Y\!\in\!H_2(Y;\Z)$ is the homology class of 
a line in the exceptional divisor,  $\be\!=\!\io_{Y*}(\be_Y)$, $a_1,a_2\!=\!0$, 
and $f_1,f_2\!:\P^1\!\lra\!X$ are parametrizations of proper transforms of two distinct lines
in~$\P^n$ passing through the center of the blowup.
In this particular case, 
$$\ov\M_{0,\f}(X,\be;J)=\ov\M_{0,\f}(Y,\be_Y;J)$$
and $\cok(D_J^{NY})$ is the zero vector bundle.
Thus, if $\bar\ell$ denotes the homology class of $f_1$ and~$f_2$,
$$\big(\bar\ell,\bar\ell\big)_{0,\be}^X=
\big((\tau_{a_j}\ka_j)_{j\in S}\big)_{0,\be}^X
=\bC_{0,\f}(Y,\be_Y)
=\,^{\pm}\big|\ov\M_{g,\f}(Y,\be_Y;J)\big|
=\big(\bar\ell\cap Y,\bar\ell\cap Y\big)_{0,\be_Y}^Y
=1;$$
see Footnote~\ref{inY_ft}.
\end{eg}

\noindent
Various special cases of Theorem~\ref{main_thm},
such as those in Examples~\ref{CY_eg}-\ref{transversetoY_eg},  
are standard in the algebraic setting and are used 
in \cite{BP}, \cite{KlP}, and~\cite{P1}, for example.
Some special cases of Theorem~\ref{main_thm} have appeared in 
the symplectic setting as well, including in \cite{LiZ}, 
\cite{McTo}, and~\cite{Taubes}.
Examples~\ref{inY_eg} and~\ref{transversetoY_eg} generalize Example~\ref{CY_eg}
in two opposite directions.
Corollary~\ref{LeP_crl} below, which applies this theorem in the setting of~\cite{LeP},
is yet another special case of Example~\ref{transversetoY_eg}.
The full statement of Theorem~\ref{main_thm} mixes the two extreme cases
of Examples~\ref{inY_eg} and~\ref{transversetoY_eg}.\\

\noindent 
The striking conclusion of~\cite{LeP} is that all GW-invariants 
of a Kahler surface~$X$ of general type localize to a canonical divisor.
The situation is particularly beautiful if $X$ admits a smooth canonical divisor~$\K_X$.
If $X$ is minimal, the GW-invariants of~$X$ in degrees other than multiples of~$\K_X$ vanish.
The GW-invariants of~$X$ in degrees~$\K_X$ and~$2\K_X$ are computed in~\cite{KLi} 
via an algebraic reformulation of~\cite{LeP}
and shown to satisfy a conjecture of~\cite{MaP}.
In the next paragraph we review the relevant statements from~\cite{LeP}.\\

\noindent
Let $(X,J_0)$ be a minimal Kahler surface of general type and $\al$
the real part of a non-zero holomorphic $(2,0)$-form such that
$Y\!\equiv\!\al^{-1}(0)$ is smooth (and reduced).
Since $X$ is minimal, $Y$ is connected.
With $\lr{\cdot,\cdot}$ denoting the Riemannian metric on~$X$, define
\begin{gather}
K_{\al}\in \Ga\big(X;\Hom_{\R}(TX,TX)\big),  \quad
R_{\al}\in \Ga\big(Y;\Hom_{\R}(TY\!\otimes\!_{\C}NY,NY)\big),
\qquad\hbox{by}\notag\\
\label{Rdfn_e}\begin{split}
\lr{v_1,K_{\al}v_2}&=\al(v_1,v_2)~~\forall\,v_1,v_2\!\in\!T_xX,\,x\!\in\!X;\\
R_{\al}(v_1,v_2)&= J_0\big\{\na_{v_2}K_{\al}\big\}(v_1)+T_xY
~~\forall\,v_1\!\in\!T_xY,\,v_2\!\in\!T_xX,\,x\!\in\!X.
\end{split}
\end{gather}
By \cite[Lemmas 2.1,8.2]{LeP}, $R_{\al}$ is well-defined.
The almost complex structure~$J_{\al}$ on $X$ described in \cite[Section~2]{LeP}
agrees with~$J_0$ along the smooth complex curve~$Y$.
By \cite[Lemma~2.3]{LeP}, every non-constant $J_{\al}$-holomorphic map 
$u\!:\Si_u\!\lra\!X$ is in fact a $J_0$-holomorphic map to~$Y$
and so lies in the homology class $dY$ for some $d\!\in\!\Z^+$.
By \cite[Section~8]{LeP}, the operator on the normal bundle~$NY$ of~$Y$
induced by the linearization of the $\dbar_{J_{\al}}$-operator for maps to~$X$
at such a map~$u$ is given~by 
\BE{DSurf_e} D_{J_{\al};u}^{NY}\!=\!\dbar_{u^*NY}+R_{\al}(du,\cdot)\!:
L^p_1(\Si_u;u^*NY)\lra 
L^p\big(\Si_u;T^*\Si_u^{0,1}\!\otimes_{\C}\!u^*NY\big),\EE
where $\dbar_{u^*NY}$ is the $\dbar$-operator in the holomorphic bundle
$u^*(NY,J_0)\!\lra\!\Si_u$.
By \cite[Proposition~8.6]{LeP}, $D_{J_{\al};u}^{NY}$ is injective.
In light of Theorem~\ref{main_thm},
Corollary~\ref{LeP_crl} below is thus simply a re-formulation
of the main conclusion of~\cite{LeP}.

\begin{crl}\label{LeP_crl}
Suppose $(X,J_0)$ is a minimal Kahler surface of general type, 
$\al$ is the real part of a non-zero holomorphic $(2,0)$-form such that
$Y\!\equiv\!\al^{-1}(0)$ is smooth, $g\!\in\!\bar\Z^+$, 
$d\!\in\!\Z^+$, $S$ is a finite set,
$S_2\!\subset\!S$, $a_j\!\in\!\bar\Z^+$ for each $j\!\in\!S$,  
and $\ka_j\!\in\!H_2(X;\Z)$ for each $j\!\in\!S_2$.
If $R_{\al}$ is defined by~\e_ref{Rdfn_e}, then the cokernels
of the operators~\e_ref{DSurf_e} form a natural oriented vector orbi-bundle
$$\cok\big(D_{\al}^{NY}\big) \lra \ov\M_{g,S}(Y,dY)$$
and
\begin{equation*}\begin{split}
&\big((\tau_{a_j}\ka_j)_{j\in S_2},(\tau_{a_j}1)_{j\in S-S_2}\big)_{g,d\K_X}^X\\
&\qquad
=\bigg(\prod_{j\in S_2}\lr{c_1(T^*X),\ka_j}\bigg)
\bigg\lan e\big(\cok(D_{\al}^{NY})\big)\!\prod_{j\in S_2}\!\!\!\big(\ev_j^*PD_Y(pt)\big)
\prod_{j\in S}\!\psi_j^{a_j},
\big[\ov\M_{g,S}(Y,dY)\big]^{vir}\bigg\ran.
\end{split}\end{equation*}
\end{crl}

\subsection{The Fano case of the Gopakumar-Vafa prediction}
\label{FanoGV_subs}

\noindent
GW-invariants are generally not integers.
On the other hand, at least in the case of projective $3$-folds (symplectic $6$-manifolds),
certain combinations of them are believed to be integers.
Ideally these combinations would be precisely counts of curves
of fixed genus and degree and passing through appropriate constraints.
A projective $3$-fold~$X$ is never ideal in this sense, but one might hope 
that $X$ becomes ideal if its Kahler complex structure is replaced
with a generic almost complex one.
We show that this is indeed the case in the ``Fano" case.\\

\noindent
If $(X,\om)$ is a compact symplectic manifold, $g\!\in\!\bar\Z^+$, $S$ is a finite set,
$\be\!\in\!H_2(X;\Z)$, and $J$ is an $\om$-tame almost complex structure on~$X$, let
$$\M_{g,S}^*(X,\be;J)\subset\ov\M_{g,S}(X,\be;J)$$
denote the subspace consisting of \sf{simple} maps, i.e.~$J$-holomorphic maps
$u\!:\Si_u\!\lra\!X$ such that $\Si_u$ is a smooth (connected) Riemann surface and 
$u^{-1}(u(z))\!=\!\{z\}$ and $d_zu\!\neq\!0$ for some $z\!\in\!\Si_u$.
These conditions imply that $u$ does not factor through a $d$-fold cover
$\Si_u\!\lra\!\Si$, with $d\!>\!1$; see~\cite[Section~2.5]{McS}.
If $f_j\!: M_j\!\lra\!X$, $j\!\in\!S$, are smooth maps from compact oriented manifolds
of even dimensions, let
$$\M_{g,\f}^*(X,\be;J)= \ov\M_{g,\f}(X,\be;J)\cap
\bigg(\M_{g,S}^*(X,\be;J)\times \prod_{j\in S}M_j\bigg),$$
with $\ov\M_{g,\f}(X,\be;J)$ defined by~\e_ref{Mfdfn_e}.
If $\M_{g,\f}^*(X,\be;J)$ is a finite set consisting of regular pairs $([u],(w_j)_{j\in S})$,
we will denote its signed cardinality by $E_{g,\be}^X(J,\f)$.\\

\noindent
If the (real) dimension of $X$ is $6$, the expected dimension of the moduli space 
$\ov\M_{g,S}(X,\be;J)$ is independent of the genus~$g$; see~\e_ref{virdim_e}.
Thus, one can mix curve counts of different genera passing through 
the same constraints.
Furthermore, if $\be\!\in\!H_2(X;\Z)$ and $\lr{c_1(TX),\be}\!<\!0$,
all degree~$\be$ GW-invariants are zero, since 
the moduli space of unmarked maps has negative expected dimension.
This leaves the ``Calabi-Yau" case, $\lr{c_1(TX),\be}\!=\!0$, and
the ``Fano" case, $\lr{c_1(TX),\be}\!>\!0$.
If $g,h\!\in\!\bar\Z^+$, define $C_{h,\be}^X(g)\!\in\!\Q$ by
\BE{Cdfn_e} \sum_{g=0}^{\i}C_{h,\be}^X(g)t^{2g}
=\bigg(\frac{\sin(t/2)}{t/2}\bigg)^{2h-2+\lr{c_1(TX),\be}}\,.\EE

\begin{thm}\label{FanoGV_thm}
Suppose $(X,\om)$ is a compact symplectic $6$-fold, $\be\!\in\!H_2(X;\Z)$,
$g\!\in\!\bar\Z^+$, $S$ is a finite set, 
and $\ka_j\!\in\!H_*(X;\Z)$ for $j\!\in\!S$ are
such that~\e_ref{dimcond_e} is satisfied with $a_j\!=\!0$.
If $\lr{c_1(TX),\be}\!>\!0$,
\begin{enumerate}[label=(\arabic*)]
\item there exists a dense open subset $\cJ_{\reg}(g,\be)$
of the space of smooth $\om$-tame almost complex structures on~$X$
such that for all $h\!\le\!g$:
\begin{itemize}
\item the moduli space $\M_{h,S}^*(X,\be;J)$ consists of regular maps;
\item for a generic choice of pseudocycle representatives\footnotemark 
$f_j\!:M_j\!\lra\!X$ for $\ka_j$, $\M_{h,\f}^*(X,\be;J)$ is a finite set of
regular pairs $([u],(w_j)_{j\in S})$ such that $u$ is an embedding;
\end{itemize}
\item the numbers $E_{h,\be}^X(\f,J)$, with $h\!\le\!g$, are independent of the choice of 
$J\!\in\!\cJ_{\reg}(g,\be)$ and~$f_j$ and can thus be denoted 
$E_{h,\be}^X((\ka_j)_{j\in S})$;
\item if $C_{g,\be}^X(h)$ is defined by~\e_ref{Cdfn_e},
\BE{FanoGV_e} \big((\ka_j)_{j\in S}\big)_{g,\be}^X=
\sum_{h=0}^{h=g} C_{h,\be}^X(g\!-\!h) E_{h,\be}^X\big((\ka_j)_{j\in S}\big).\EE\\
\end{enumerate}
\end{thm}
\footnotetext{After replacing $\ka_j$ by a multiple, $M_j$ can be taken to
be a smooth compact manifold.}

\noindent
For $g\!=\!0,1$, \e_ref{FanoGV_e} gives
\BE{lowgenusGV_e}\begin{split}
\big((\ka_j)_{j\in S}\big)_{0,\be}^X&=E_{0,\be}^X((\ka_j)_{j\in S}\big), \\
\big((\ka_j)_{j\in S}\big)_{1,\be}^X&=E_{1,\be}^X((\ka_j)_{j\in S}\big)
+\frac{2-\lr{c_1(TX),\be}}{24}E_{0,\be}^X((\ka_j)_{j\in S}\big).
\end{split}\EE
The first identity expresses the well-known fact that the genus $0$
GW-invariants of a Fano manifold are enumerative.
The second identity in~\e_ref{lowgenusGV_e} is the $n\!=\!3$ case of the relation 
between  the standard genus~$1$ GW-invariants and the reduced genus~$1$ GW-invariants
constructed in~\cite{g1comp2} for all symplectic manifolds.\\

\noindent
By the proof of \cite[Theorem~3.1.5]{McS}, 
for a generic almost complex structure $J$ on $X$ all moduli spaces 
$\M_{h,\eset}^*(X,\be';J)$ are smooth and of the expected dimension, $2\lr{c_1(TX),\be'}$.
In particular,
\BE{emptyspace_e} \lr{c_1(TX),\be'}<0 \qquad\Lra\qquad
\M_{h,S}^*(X,\be';J),\ov\M_{h,S}(X,\be';J)=\eset .\EE
By a similar argument, for a generic $J$ on $X$ the evaluation maps
$$\ev_1,\ev_2\!:  \M_{g,\{1,2\}}^*(X,\be;J)\lra X$$
are transverse, while the bundle section
$$\M_{g,\{1\}}^*(X,\be;J)\lra L_1^*\!\otimes\!\ev_1^*TX, \qquad
[u]\lra d_{z_1(u)}u\,,$$
where $L_1\!\lra\!\M_{g,\{1\}}^*(X,\be;J)$ is the universal tangent line bundle
at the marked point and $z_1(u)\!\in\!\Si_u$ is the marked point of~$u$,
is transverse to the zero set.
Thus, 
$$\M_{g,S}^{sing}(X,\be;J)\equiv
\big\{[u]\!\in\!\M_{g,S}^*(X,\be;J)\!:~
u~\textnormal{is not an embedding}\big\}$$
is the image of a smooth map from a smooth manifold of (real) dimension 
two less than the dimension of $\M_{g,S}^*(X,\be;J)$.
It follows that for a generic choice of pseudocycle representatives $f_j\!:M_j\!\lra\!X$
for~$\ka_j$, $\M_{g,\f}^*(X,\be;J)$ is a $0$-dimensional oriented sub-manifold of
$$\big(\M_{g,S}^*(X,\be;J)-\M_{g,S}^{sing}(X,\be;J)\big)
\times\prod_{j\in S}\!M_j.$$\\

\noindent
We next show that $\M_{g,\f}^*(X,\be;J)$ is a finite set.
If not, there is a sequence $([u_r],(w_{r,j})_{j\in S})$ in $\M_{g,\f}^*(X,\be;J)$
converging to some
$$\big([u],(w_j)_{j\in S}\big)\in 
\ov\M_{g,\f}(X,\be;J)- \M_{g,\f}^*(X,\be;J).$$
The image of $u$ is a connected $J$-holomorphic curve in $X$ of genus $h\!\le\!g$ 
with $k\!\ge\!1$ irreducible components of degrees $\be_1,\ldots,\be_k\!\in\!H_2(X;\Z)$
such that 
$$ d_1\be_1+\ldots+d_k\be_k=\be \qquad\hbox{for some}\quad
d_1,\ldots,d_k\in\Z^+.$$
By~\e_ref{emptyspace_e}, $\lr{c_1(TX),\be_i}\!\ge\!0$ for all $i\!=\!1,\ldots,k$.
Thus,
$$\sum_{i=1}^{i=k}\lr{c_1(TX),\be_i}\le \lr{c_1(TX),\be}. $$
The dimension-counting argument  of~\cite[Section~6.6]{McS} then shows that 
$k\!=\!1$ and $d_1\!=\!1$.
It then follows that the image of~$u$ is an irreducible $J$-holomorphic curve
of degree~$\be$ and genus~$h\!<\!g$ that meets each of the maps $f_j$ with $j\!\in\!S$.\\

\noindent
While degree~$\be$ genus~$h\!<\!g$ $J$-holomorphic curves meeting the maps $f_j$ 
can certainly exist for a generic $J$, they cannot be limits of other degree~$\be$
curves meeting the maps~$f_j$ by the $\nu_r\!=\!0$ case of Proposition~\ref{horreg_prp}
for the following reason.
If
$$\big([u],(w_j)_{j\in S}\big)\in\ov\M_{g,\f}(X,\be;J)
-\M_{g,\f}^*(X,\be;J),$$
the domain of $u$ consists of two or more irreducible components.
Furthermore, by the previous paragraph, the restriction of~$u$ to all components,
except for one, is constant;
let $u_{\eff}$ denote the \sf{effective part} of~$u$, 
i.e.~the non-constant restriction.
The domain $\Si_{u_{\eff}}$ of~$u_{\eff}$ is a smooth curve of genus~$h\!<\!g$ 
with distinct points
$(z_j(u_{\eff}))_{j\in S}$ that are mapped to $(\ev_j(u))_{j\in S}$ by~$u_{\eff}$.
Thus,
$$\big([u_{\eff}],(w_j)_{j\in S}\big)\in\M_{h,\f}^*(X,\be;J);$$
by the previous paragraph, $u_{\eff}$ is an embedding onto a smooth $J$-holomorphic curve~$Y$
of genus~$h$ degree~$\be$ meeting the maps~$f_j$.
This implies that removing a node from $\Si_{u_{\eff}}$ 
disconnects~$\Si_u$.\footnote{This observation implies that the homomorphism
\e_ref{kerrestr_e} is surjective.}
Since the total evaluation map
$$\Bev\!\equiv\!\prod_{j\in S}\ev_j: \M_{h,S}^*(X,\be;J)\lra X^S$$
is transverse to $\f$, 
\BE{Ycount_e}\ker(D_{J;u_{\eff}}^{NY})\lra \bigoplus_{j\in S}N_{f_j(w_j)}^{f_j}Y, \qquad
\xi\lra \big(\xi(z_j(u_{\eff}))+T_{f_j(w_j)}Y+(\Im\,d_{w_j}f_j)\big)_{j\in S}\,,\EE
is surjective; see Section~\ref{AbsGW_subs} for the notation.
Since $u_{\eff}$ is a regular map,
\begin{equation*}\begin{split}
\dim \ker(D_{J;u_{\eff}}^{NY})
=\ind\big(D_{J;u_{\eff}}^{NY}\big)
&=2\big(\lr{c_1(NY),Y}+2(1\!-\!h)\big)=2\lr{c_1(TX),\be}\\
&=\sum_{j\in S}(4\!-\!\dim\,M_j)
\le\sum_{j\in S}\dim\,N^{f_j}_{f_j(w)}Y;\end{split}\end{equation*}
the second-to-last equality holds by~\e_ref{dimcond_e}.
Thus, the homomorphism in~\e_ref{Ycount_e} is an isomorphism.
On the other hand, $D_{J;u}^{NY}$ is the restriction 
of the operator $\bigoplus_i D_{J;u_i}^{NY}$ to 
$$L^p_1(\Si_u;u^*NY)\subset\bigoplus_i L^p_1\big(\Si_{u;i};u_i^*NY\big),$$
where $\{\Si_{u;i}\}$ are the irreducible components of~$\Si$ and $u_i\!=\!u|_{\Si_{u;i}}$.
If $u_i$ is a constant map, then $D_{J;u_i}^{NY}$ is the usual $\dbar$-operator 
on the space of functions on $\Si_{u_i}$ with values in $N_{u_i(\Si_{u;i})}Y\!\approx\!\C^2$.
Since $\Si_u$ is a connected nodal Riemann containing~$\Si_{u_{\eff}}$ as a component,
$u|_{\Si_{\eff}}\!=\!u_{\eff}$, and $u$ is constant on each of the irreducible
components of  $\Si_u\!-\!\Si_{u_{\eff}}$, it follows that the projection homomorphism
\BE{kerrestr_e} \ker D_{J;u}^{NY}\lra\ker D_{J;u_{\eff}}^{NY}, \qquad
\xi\lra\xi|_{\Si_{u_{\eff}}},\EE
is an isomorphism.
Thus, the homomorphism
$$\ker(D_{J;u}^{NY})\lra \bigoplus_{j\in S}N_{f_j(w_j)}^{f_j}Y, \qquad
\xi\lra \big(\xi(z_j(u))+T_{f_j(w_j)}Y+(\Im\,d_{w_j}f_j)\big)_{j\in S}\,,$$
is an isomorphism, since the homomorphism \e_ref{Ycount_e} is. 
Therefore, by Proposition~\ref{horreg_prp} there is no sequence in
$$\ov\M_{g,\f}(X,\be;J)-\ov\M_{g,\f}(Y,[Y];J)\supset \M_{g,\f}^*(X,\be;J)$$
converging to $([u],(w_j)_{j\in S})$.\\

\noindent
We have thus shown that $\M_{g,\f}^*(X,\be;J)$ is a compact oriented $0$-dimensional 
manifold and its signed cardinality $E_{g,\be}^X(\f,J)$ is well-defined.
The independence of $E_{g,\be}^X(\f,J)$ of the choices of $J$ and $f_j$
follows from \e_ref{FanoGV_e}, with $E_{h,\be}^X((\ka_j)_{j\in S})$
replaced by $E_{h,\be}^X(\f,J)$.
In turn, this identity follows from Theorem~\ref{main_thm} and 
the proof of \cite[Theorem~3]{P1}.
Let $Y$ be a degree~$\be$ $J$-holomorphic curve of genus~$h\!\le\!g$
meeting each~$f_j$.
By the above, the assumptions of Theorem~\ref{main_thm} are satisfied.
By definition (see Section~\ref{CR_subs2}), 
the orbi-bundle $\cok(D_J^{NY})$ is dual to the bundle $\ker((D_J^{NY})^*)$ of 
kernels of the dual operators  $(D_J^{NY})^*$.
For each
$$\big([u],(w_j)_{j\in S}\big)\in\ov\M_{g,\f}(Y,[Y];J)
\subset\ov\M_{g,\f}(X,\be;J),$$
the operator $(D_{J;u}^{NY})^*$ is the natural extension of the operator 
$\bigoplus_i (D_{J;u_i}^{NY})^*$
to $(1,0)$-forms on $\Si_u$ with poles at the nodes such that the residues
at each node sum up to~$0$.
Since $(D_{J;u_{\eff}}^{NY})^*$ is injective by the regularity of~$u_{\eff}$,
the projection 
$$\eta\lra \bigoplus_{\Si_{u;i}\neq\Si_{u_{\eff}}}\eta|_{\Si_{u;i}}$$
to the contracted components is injective. 
Since $(D_{J;u_i}^{NY})^*\!=\!\dbar^*$ if $u_i$ is constant,
the image of this homomorphism is determined by~$\Si_u$ 
and is independent of $D_{J;u_{\eff}}^{NY}$ 
(as long as $D_{J;u_{\eff}}^{NY}$ is surjective).
Thus, $\cok(D_J^{NY})$ is isomorphic to the restriction to $\ov\M_{g,\f}(Y,[Y];J)$
of the obstruction bundle in \cite[Section~3]{P1},
i.e.~the bundle of cokernels of the operators~$D_{J;u}^{NY}$ as above,
but for a holomorphic vector bundle~$NY$.
Thus,
\BE{P2_e}\begin{split}
\bC_{g,\f}(Y,\be_Y)
&=\bigg\lan e\big(\cok(D_J^{NY})\big),
\big[\ov\M_{g,\f}(Y,[Y];J)\big]^{vir}\bigg\ran\\
&=C_{h,\be}^X(g\!-\!h)\,\sgn\big([u_{\eff}],(w_j)_{j\in S}\big)
\end{split}\EE
by \e_ref{mainthm_e2} and \cite[Theorem~3]{P1}.
Since
$$\ov\M_{g,\f}(X,\be;J)
=\bigsqcup_{h=0}^{h=g}\bigsqcup_{([u],(w_j)_{j\in S})\in\M_{h,\f}^*(X,\be;J)}
\!\!\!\!\!\!  \ov\M_{g,\f}(\Im\,u,[\Im\, u];J),$$
the identity~\e_ref{FanoGV_e} follows from~\e_ref{P2_e}.\\

\noindent
Theorem~\ref{FanoGV_thm} confirms (a stronger version of) 
the Fano case of \cite[Conjecture~2(i)]{P2},
i.e.~that the numbers $E_{h,\be}^X((\ka_j)_{j\in S})$
{\it defined from GW-invariants by~\e_ref{FanoGV_e}} are integers.
The Calabi-Yau case is fundamentally more difficult as it involves 
multiple covers of curves.\footnote{Theorem~\ref{FanoGV_thm} and its proof also apply 
to the cases when $\lr{c_1(TX),\be}\!=\!0$, but $\be$ is not a non-trivial integer multiple
of another element of $H_2(X;\Z)$.}
On the other hand, it might be possible to approach \cite[Conjecture~2(ii)]{P2},
i.e.~that $E_{h,\be}^X((\ka_j)_{j\in S})\!=\!0$ for a fixed $\be$ and all sufficiently
large~$g$ if $X$ is projective, by studying possible limits of $J_t$-holomorphic 
curves with $J_t\!\in\!\cJ_{\reg}(g,\be)$ as $J_t$ approaches the standard complex
structure on $X\!\subset\!\P^n$ and using the Castelnuovo bound \cite[p116]{ACGH}.\\

\noindent
An algebro-geometric approach to Theorem~\ref{FanoGV_thm} has recently been proposed 
in~\cite{KKO}, at least in the usual, more narrow, meaning of {\it Fano} in algebraic
geometry.
The stable-map style invariants of smooth projective varieties defined in~\cite{KKO}
are a priori integers in the case of Fano varieties, 
just like the numbers $E_{h,\be}^X((\ka_j)_{j\in S})$.
In addition, in this Fano case, they are non-negative integers and satisfy
the vanishing prediction of \cite[Conjecture~2(ii)]{P2}.
However, it remains to be shown that they are related to the GW-invariants 
in the required way, i.e.~as in~\e_ref{FanoGV_e}.

\section{Analytic Preliminaries}
\label{analysis_sec}

\noindent
In this section, we collect a number of background statements concerning 
solutions of perturbed Cauchy-Riemann equations.
For the rest of the paper, fix a real number $p\!>\!2$.
If $\Si$ is a  $2$-dimensional manifold, this condition implies
that any $L^p_1$-map $\Si\!\lra\!\R$ is continuous and 
in particular has a well-defined value at each point.

\subsection{Nodal Riemann surfaces}
\label{NRS_subs}

\noindent
Let $(E,\fI)\!\lra\!\Si$ be an $L^p_1$-complex vector bundle over 
a smooth Riemann surface, i.e.~a one-dimensional complex manifold.
If $z\!\in\!\Si$ and
$$A_z\in\Hom_{\R}(E_z,T_z^*\Si^{0,1}\!\otimes_{\C}\!E_z),$$
we define
\begin{gather*}
A_z^*\in\Hom_{\R}(T_z^*\Si^{1,0}\!\otimes_{\C}\!E_z^*,T_z^*\Si^{1,1}\!\otimes_{\C}\!E_z^*)
\qquad\hbox{by}\\
\Re\big(v\wedge (A_z^*w)\big)= \Re\big((A_zv)\wedge w\big)
\in \La_{\R}^2(T_z^*\Si)
\qquad\forall\, v\!\in\!E_z,\,w\!\in\!T_z^*\Si^{1,0}\!\otimes_{\C}\!E_z^*\,.
\end{gather*}
Since $\La_{\R}^2(T_z^*\Si)$ is one-dimensional, $A_z^*$ is well-defined.
If
$$A\in L^p\big(\Si;\Hom_{\R}(E,T^*\Si^{0,1}\!\otimes_{\C}\!E)\big),$$
this construction gives rise to an element 
\begin{gather}
A^*\in L^p\big(\Si;\Hom_{\R}(T^*\Si^{1,0}\!\otimes_{\C}\!E^*,
T^*\Si^{1,1}\!\otimes_{\C}\!E^*)\big) \qquad\hbox{s.t.} \notag\\
\label{Aadj_e}
\bllrr{\xi,A^*\eta}\equiv
\Re\Big(\int_{\Si}\xi\wedge (A^*\eta)\Big)
=\Re\Big(\int_{\Si}(A\xi)\wedge\eta\Big)
\equiv \bllrr{A\xi,\eta}
\end{gather}
for all $\xi\!\in\!L^p_1(\Si;E)$ and
$\eta\!\in\!L^p_1(\Si;T^*\Si^{1,0}\!\otimes\!E^*)$.\\

\noindent
Let $E\!\lra\!\Si$ be as above.
If  $S$ is a finite subset of~$\Si$, denote~by 
$$L_k^p\big(\Si;E(S)\big)\subset L^p_{k,loc}(\Si\!-\!S;E)$$
the subspace of sections $\eta$ of $E$ such that for every 
$z_0\!\in\!S$ there exist a neighborhood~$U$ of~$z_0$ in~$\Si$ and
a coordinate $w\!:U\!\lra\!\C$ such that 
$$w(z_0)=0 \qquad\hbox{and}\qquad
w\cdot\eta|_U\in L_k^p(U;E).$$ 
If $k\!\ge\!1$, 
an element~$\eta$ of $L_k^p(\Si;T^*\Si^{1,0}\!\otimes_{\C}\!E(S))$ 
has a well-defined residue at $z_0\!\in\!S$ given~by
$$\Res_{z=z_0}\eta=\xi(z_0)\in E_{z_0} \qquad\hbox{if}\qquad
\eta(z)=\frac{dw}{w(z)}\otimes\xi(z)~~\forall~z\!\in\!U,~
\xi\in L_1^p(U;E).\footnotemark$$
If $\vr$ is a function assigning to each element $z_0\!\in\!S$ 
a real subspace $E_{z_0}'\!\subset\!E_{z_0}$, let
$$L_1^p\big(\Si;T^*\Si^{1,0}\!\otimes_{\C}\!E(\vr)\big)=\big\{
\eta\!\in\!L_1^p\big(\Si;E(S)\big)\!: \,
\Res_{z=z_0}\eta\!\in\!E_{z_0}'~\forall\,z_0\!\in\!S\big\}.$$\\
\footnotetext{If $\eta\in L_k^p(\Si;T^*\Si^{1,0}\!\otimes_{\C}\!E(S\!-\!z_0))$,
then $\Res_{z=z_0}\eta=0$. The converse is not true; for example, the residue of
$\eta\!=\!\bar{z}\,dz/z$ is zero at $z\!=\!0$, but $\eta$ is not even continuous 
at $z\!=\!0$.
On the other hand, the converse is true if $\eta$ lies in the kernel of
a generalized CR-operator as in Section~\ref{CR_subs}.}

\noindent
By a \sf{Riemann surface} $\Si$ we will mean a compact complex one-dimensional
manifold with pairs of distinct points identified.
In other words, 
\BE{RSdfn_e}\Si=\ti\Si\big/\sim, \qquad\hbox{where}\quad
x_i^{(1)}\sim x_i^{(2)}~~i=1,\ldots,m,\EE
for some smooth compact Riemann surface $\ti\Si$ and 
distinct points~$x_i^{(1)},x_i^{(2)}\!\in\!\ti\Si$.
The quotient map
$$\si\!: \ti\Si\lra\Si$$
is determined by $\Si$ up to an isomorphism.
We will denote by 
$$\Si_{\sing}\equiv\big\{\si(x_i^{(1)})\!:\,i\!=\!1,\ldots,m\big\}
\subset\Si
\qquad\hbox{and}\qquad
\ti\Si_{\sing}\equiv\big\{x_i^{(1)},x_i^{(2)}\!:\,i\!=\!1,\ldots,m\big\}
\subset\ti\Si$$
the subset of \sf{singular points} of $\Si$ and its preimage under $\si$,
respectively.
Let $\Si^*\!\subset\!\Si$ be the subspace of \sf{smooth points},
i.e.~the complement of~$\Si_{\sing}$.\\

\noindent
If $Y$ is a smooth manifold and~$\Si$ is a Riemann surface as above, 
an \sf{$L^p_1$-map $u\!:\Si\!\lra\!Y$} is an $L^p_1$-map
$$\ti{u}\!:\ti\Si\lra Y \qquad\hbox{s.t.}\quad 
\ti{u}\big(x_i^{(1)}\big)=\ti{u}\big(x_i^{(2)}\big)~~\forall\,i=1,\ldots,m.$$
By a \sf{vector bundle} $E\!\lra\!\Si$, 
we will mean a topological complex vector bundle such that
$\si^*E\!\lra\!\ti\Si$ is an $L^p_1$-complex vector bundle.
Let
\begin{equation*}\begin{split}
L_1^p(\Si;E)
&=\big\{\xi\!\in\!L_1^p(\ti\Si;\si^*E\big)\!:\,
\xi(x_i^{(1)})\!=\!\xi(x_i^{(2)})~\forall\,i\!=\!1,\ldots,m\big\};\\
L^p\big(\Si;T^*\Si^{0,1}\!\otimes\!_{\C}E\big)
&=L^p\big(\ti\Si;T^*\ti\Si^{0,1}\!\otimes\!_{\C}\si^*E\big).
\end{split}\end{equation*}
If $S$ is a finite subset of~$\Si^*$, let $\ti{S}\!=\!\si^{-1}(S)$ and define
\begin{equation}\label{dualsheafdfn_e}\begin{split}
L_1^p\big(\Si;\K_{\Si}\!\otimes_{\C}\!E(S)\big)
&=\Big\{\eta\!\in\!
L_1^p\big(\ti\Si;T^*\ti\Si^{1,0}\!\otimes_{\C}\!
\si^*E(\ti{S}\!\cup\!\ti\Si_{\sing})\big)\!:
\\ &\qquad\qquad\qquad
\sum_{\ti{z}_0\in\si^{-1}(z_0)}\!\!\!\!\!\!\Res_{z=\ti{z}_0}\eta(\ti{z}_0)\!=\!0~~\forall\,
z_0\!\in\!\Si_{\sing}\Big\},\\
L^p\big(\Si;T^*\Si^{0,1}\!\otimes_{\C}\!\K_{\Si}\!\otimes\!_{\C}E(S)\big)
&=L^p\big(\ti\Si;T^*\ti\Si^{0,1}\!\otimes_{\C}\!T^*\ti\Si^{1,0}\!\otimes_{\C}\!
\si^*E(\ti{S}\!\cup\!\ti\Si_{\sing})\big).
\end{split}\end{equation}
If $\vr$ is a function assigning to each element $z_0\!\in\!S$ 
a real subspace $E_{z_0}'\!\subset\!E_{z_0}$, let
\BE{dualsheafdfn_e2} L_1^p\big(\Si;\K_{\Si}\!\otimes_{\C}\!E(\vr)\big)=\big\{
\eta\!\in\!L_1^p\big(\Si;\K_{\Si}\!\otimes_{\C}\!E(S)\big)\!: \,
\Res_{z=\si^{-1}(z_0)}\eta\!\in\!E_{z_0}'~\forall\,z_0\!\in\!S\big\}.\EE
Similarly, we define
\begin{equation*}\begin{split}
L_1^p\big(\Si;E(-S)\big)
&=\big\{\xi\!\in\!L_1^p\big(\Si;E)\!:\,\xi(z_0)\!=\!0~\forall\,z_0\!\in\!S\big\},\\
L_1^p\big(\Si;E^*(-\vr)\big)
&=\big\{\xi\!\in\!L_1^p\big(\Si;E^*)\!:\,\xi(z_0)\!\in\!\Ann(E_{z_0}')
~\forall\,z_0\!\in\!S\big\},
\end{split}\end{equation*}
where $\Ann(E_{z_0}')\!\subset\!\Hom_{\R}(E_{z_0},\R)$ 
is the annihilator of $E_{z_0}'\!\subset\!E_{z_0}$.
The real pairings in~\e_ref{Aadj_e} extend to pairings
\begin{equation*}\begin{split}
L_1^p\big(\Si;E\big)\!\otimes\!
L^p\big(\Si;T^*\Si^{0,1}\!\otimes_{\C}\!\K_{\Si}\!\otimes_{\C}\!E^*(S)\big)
&\lra\R,\\
L^p\big(\Si;T^*\Si^{0,1}\!\otimes_{\C}\!E\big)\!\otimes\!
L_1^p\big(\Si;\K_{\Si}\!\otimes_{\C}\!E^*(S)\big) 
&\lra\R.
\end{split}\end{equation*}
Furthermore, the equality in \e_ref{Aadj_e} holds for all
$\eta\in L_1^p\big(\Si;\K_{\Si}\!\otimes_{\C}\!E^*(S)\big)$.

\subsection{Generalized Cauchy-Riemann operators}
\label{CR_subs}

\begin{dfn}\label{CR_dfn}
Let $(Y,J)$ be an almost complex manifold and 
$(N,\fI)\!\lra\!(Y,J)$ a smooth vector bundle.
\begin{enumerate}[label=(\arabic*)]
\item A \sf{$\dbar$-operator on $(N,\fI)$} is a $\C$-linear map
$$\dbar\!: \Ga(Y;N)\lra \Ga^{0,1}(Y;N)\equiv
\Ga\big(Y;T^*Y^{0,1}\!\otimes\!_{\C}N\big)$$
such that 
$$\dbar\big(f\xi)=(\dbar{f})\!\otimes\!\xi+f(\dbar\xi)
\qquad\forall~f\!\in\!C^{\i}(Y),~\xi\!\in\!\Ga(Y;N).$$
\item A \sf{smooth generalized Cauchy-Riemann operator} (or \sf{smooth CR-operator}) 
on $(N,\fI)$ is a differential operator of the form 
\BE{CRsplit_e0} D=\dbar\!+\!A\!: \Ga(Y;N)\lra \Ga^{0,1}(Y;N),\EE
where $\dbar$ is a $\dbar$-operator on~$(N,\fI)$ and
$$ A\in \Ga\big(Y;\Hom_{\R}(N,T^*Y^{0,1}\!\otimes_{\C}\!N)\big).$$\\
\end{enumerate}
\end{dfn}

\noindent
If $\na$ is an affine connection in $(N,\fI)$, the operator
\BE{ClinCR_e0} \Ga(Y;N)\lra \Ga^{0,1}(Y;N),
\qquad \xi\lra \frac{1}{2}\big(\na\xi+\fI\na\xi\circ J\big),\EE
is a $\dbar$-operator on~$(N,\fI)$.
Furthermore, any $\C$-linear CR-operator on $(N,\fI)$
is a $\dbar$-operator, and 
any $\dbar$-operator on~$(N,\fI)$ is of the form~\e_ref{ClinCR_e0}
for some (not unique) connection~$\na$ in~$(N,\fI)$.
In particular, $A$ in the decomposition~\e_ref{CRsplit_e0} can be assumed
to be $\C$-anti-linear.\\

\noindent
Let $\na^J$ be the $J$-linear connection in $TY$ obtained
from a Levi-Civita connection~$\na$ on~$Y$ 
and $A_Y(\cdot,\cdot)$  the \sf{Nijenhuis tensor} of~$J$: 
\begin{equation*}\begin{split}
\na^J_{\xi_1}\xi_2&=\frac{1}{2}\Big(\na_{\xi_1}\xi_2-J\na_{\xi_1}(J\xi_2)\Big)\\
A_Y(\xi_1,\xi_2)&=\frac{1}{4}\Big([\xi_1,\xi_2]+J[\xi_1,J\xi_2]+J[J\xi_1,\xi_2]
-[J\xi_1,J\xi_2]\Big)
\end{split}
\qquad\forall\,\xi_1,\xi_2\!\in\!\Ga(Y;TY).
\end{equation*}
We identify $A_Y$ with the element
$$A_Y\in \Ga\big(Y;\Hom_{\R}(TY,T^*Y^{0,1}\!\otimes_{\C}\!TY)\big),
\qquad v\lra A_Y(\cdot,v).$$
Then,
$$\dbar_Y\!\equiv \frac{1}{2}\Big(\na^J\xi+J\na^J\circ J\Big), \,
D_Y\!\equiv\! \dbar_Y\!+\!A_Y\!:
\Ga(Y;TY)\lra \Ga^{0,1}(Y;TY)$$
are a $\dbar$-operator on $TY$ and a smooth CR-operator on~$TY$, respectively.

\begin{dfn}\label{CR_dfn2}
Let $(E,\fI)$ be an $L^p_1$ complex vector bundle over a Riemann surface~$(\Si,\fJ)$.
\begin{enumerate}[label=(\arabic*)]
\item A \sf{$\dbar$-operator on $(E,\fI)$} is a $\C$-linear map
$$\dbar\!: L^p_1(\Si;E)\lra L^p\big(\Si;T^*\Si^{0,1}\!\otimes\!_{\C}E\big)$$
such that 
$$\dbar\big(f\xi)=(\dbar{f})\!\otimes\!\xi+f(\dbar\xi)
\qquad\forall~f\!\in\!C^{\i}(\Si),~\xi\!\in\!\Ga(\Si;E).$$
\item A \sf{generalized Cauchy-Riemann operator} (or \sf{CR-operator}) 
on $(E,\fI)$ is a differential operator of the form 
\BE{CRdfn_e}D=\dbar\!+\!A\!: L^p_1(\Si;E)\lra L^p\big(\Si;T^*\Si^{0,1}\!\otimes\!_{\C}E\big),
\EE
where $\dbar$ is a $\dbar$-operator on~$(E,\fI)$ and
\BE{Adfn_e} A\in L^p\big(\Si;\Hom_{\R}(E,T^*\Si^{0,1}\!\otimes_{\C}\!E)\big).\EE\\
\end{enumerate}
\end{dfn}

\noindent
If $\na$ is an affine connection in $(E,\fI)$, the operator
\BE{ClinCR_e} L^p_1(\Si;E)\lra L^p(\Si;T^*\Si^{0,1}\!\otimes_{\C}\!E),
\qquad \xi\lra \frac{1}{2}\big(\na\xi+\fI\na\xi\circ\fJ\big),\EE
is the usual $\dbar$-operator for a unique holomorphic structure in~$(E,\fI)$.
Furthermore, any $\C$-linear CR-operator is of the form~\e_ref{ClinCR_e}.\\

\noindent
If $\Si$ and $N\!\lra\!Y$ are as above, an $L^p_1$-map $u\!:\Si\!\lra\!Y$
pulls back a smooth CR-operator~$D$ on~$N$ to a CR-operator~$D_u$ on
$u^*N\!\lra\!\Si$ as follows. 
Suppose $D$ is presented as in~\e_ref{CRsplit_e0} with $\C$-anti-linear $A$ 
and $\na$ is a connection in $(N,\fI)$ inducing the corresponding $\dbar$-operator.
Let $\ti{u}\!:\ti\Si\!\lra\!Y$ be the map corresponding to~$u$ as 
in Section~\ref{NRS_subs} and
$$\ti\na\!: L^p_1(\ti\Si;\ti{u}^*N)\lra 
L^p(\ti\Si;T^*\ti\Si\!\otimes\!_{\R}\ti{u}^*N)$$
the connection induced by $\na$.
Then,
$$D_{\ti{u}}=\frac{1}{2}\Big(\ti\na+\fI\ti\na\circ\fJ\Big)
+A\circ\partial_J \ti{u},
\qquad\hbox{where}\quad
\partial_J \ti{u}=\frac{1}{2}\Big(du-J d\ti{u}\circ \fJ\Big),$$
is a generalized CR-operator on $\ti{u}^*(N,\fI)$;
$D_{\ti{u}}$ is independent of the choice of $\na$ if $u$ is $(J,\fJ)$-holomorphic.\\

\noindent
Suppose $(Y,J)$ is an almost complex manifold and $D_Y$ is as above.
If $(\Si,\fJ)$ is a Riemann surface and $u\!:\Si\!\lra\!Y$ is 
a $(J,\fJ)$-holomorphic $L^p_1$-map, then $D_{J;u}\!\equiv\!u^*D_Y$
is the linearization of the $\dbar_J$-operator on the space of $L^p_1$-maps
from~$\Si$, with complex structure fixed, to~$Y$; see \cite[Section~3.1]{McS}.
If in addition, $(Y,J)$ is an almost complex submanifold of 
an almost complex manifold~$(X,J)$, then 
$$D_{J;u}\!\equiv\!D_{J;u}^Y\!\equiv\!u^*D_Y\!:
L^p_1(\Si;u^*TY)\lra L^p\big(\Si;T^*\Si^{0,1}\!\otimes\!_{\C}u^*TY\big)$$
is the restriction of 
$$D_{J;u}^X\!\equiv\!u^*D_X\!:
L^p_1(\Si;u^*TX)\lra L^p\big(\Si;T^*\Si^{0,1}\!\otimes\!_{\C}u^*TX\big).$$
Thus, $D_{J;u}^X$ induces a CR-operator
$$D_{J;u}^{NY}\!:
L^p_1(\Si;u^*NY)\lra L^p\big(\Si;T^*\Si^{0,1}\!\otimes\!_{\C}u^*NY\big),$$
where $NY\!\equiv\!TX|_Y/TY$ is the complex normal bundle of $Y$ in $X$.\\

\noindent 
The next lemma extends Serre duality from $\dbar$-operators to 
CR-operators.
If~$D$ is as in~\e_ref{CRdfn_e}, let
$$D^*=\dbar-A^*\!: L^p_1(\Si;\K_{\Si}\!\otimes\!_{\C}E^*)\lra 
L^p(\Si;T^*\Si^{0,1}\!\otimes_{\C}\K_{\Si}\!\otimes\!_{\C}E^*);$$
see \e_ref{Aadj_e} and~\e_ref{dualsheafdfn_e} for notation.
If $S\!\subset\!\Si$ is a finite subset of smooth points of $\Si$
and $\vr$ is a function assigning to $z_0\!\in\!S$ a complex subspace 
of~$E_{z_0}^*$, 
$D^*$ extends to an operator 
$$D_{\vr}^*\!: L^p_1\big(\Si;\K_{\Si}\!\otimes\!_{\C}E^*(\vr)\big)\lra 
L^p\big(\Si;T^*\Si^{0,1}\!\otimes_{\C}\K_{\Si}\!\otimes\!_{\C}E^*(S)\big);$$
see \e_ref{dualsheafdfn_e2}.
Let $D_{\vr}$ be the restriction of $D$ to the closed subspace
$L^p_1(\Si;E(-\vr))$ of~$L^p_1(\Si;E)$.

\begin{lmm}\label{serre_lmm}
Let $D$ be a CR-operator on a complex vector bundle $(E,\fI)$
over a  Riemann surface~$(\Si,\fJ)$.
If $S$ is a finite subset of smooth points of~$\Si$ and 
$\vr$ is a function assigning to $z_0\!\in\!S$ a real subspace of~$E_{z_0}^*$, 
the homomorphism
\BE{serre_e} \cok\, D_{\vr}\lra \Hom_{\R}(\ker D_{\vr}^*,\R), \qquad
\eta\lra\bllrr{\eta,\cdot},\EE
is an isomorphism.
\end{lmm}

\noindent
{\it Proof:}
If $\Si$ is smooth and $S\!=\!\eset$, this is \cite[Lemma~2.3.2]{IvSh}.
Furthermore, by the twisting construction of \cite[Lemma~2.4.1]{Sh}\footnote{This 
construction extends the usual procedure of twisting a holomorphic vector bundle
by a divisor to generalized CR-operators; 
it can be seen as a manifestation of Carleman Similarity Principle \cite[Theorem~2.2]{FHoS}.},
the elements~$z_0$ of $S$ for which $\vr(z_0)\!=\!E_{z_0}^*$ can be omitted 
from~$S$.
In the general case, the proof of \cite[Lemma~2.3.2]{IvSh} shows that 
the homomorphisms
\BE{serre_e1}\ker D_{\vr}^*\lra \Hom_{\R}(\cok\, D_{\vr},\R),\qquad 
\ker D_{\vr}\lra \Hom_{\R}(\cok\, D_{\vr}^*,\R),\EE
induced by the pairings \e_ref{Aadj_e} are well-defined and injective.
It follows that 
$$\ind\,D_{\vr}+ \ind\,D_{\vr}^*\le0$$
and equality holds if and only if the homomorphisms \e_ref{serre_e1}
are isomorphisms.
On the other hand, if $\ti{D}_{\vr}$ and $\ti{D}_{\vr}^*$ are the operators 
corresponding to~$D_{\vr}$ and $D_{\vr}^*$ over the normalization 
$\si\!:\ti\Si\!\lra\!\Si$, dropping any matching conditions at the nodes
and the other restricting conditions at the points of~$S$, then
\begin{equation*}\begin{split}
\ind\, D_{\vr}&=\ind\,\ti{D}_{\vr}-2k m-\|\vr\|, \\
\ind\, D_{\vr}^*&=\ind\,\ti{D}_{\vr}^*-2k m-2k|S|+\|\vr\|,
\end{split} \end{equation*}
where $k$ is the complex rank of $E$,  $m$ is the number of nodes in $\Si$, and
$$\|\vr\|=\sum_{z_0\in S}\dim_{\R}\vr(z_0).$$
Since the kernel and cokernel of~$\ti{D}_{\vr}^*$ are isomorphic to 
the kernel and cokernel of a CR-operator on $T^*\ti\Si\!\otimes\!\si^*E^*$
twisted by the preimages of the nodes and the elements of~$S$,
$$\ind\,\ti{D}_{\vr}^*=-\ind\,\ti{D}_{\vr}+4km+2k|S|.$$
It follows that $\ind\,D_{\vr}^*=-\ind\,D_{\vr}$ and thus the injective homomorphisms
in~\e_ref{serre_e1} are in fact isomorphisms.

\subsection{Families of nodal Riemann surfaces}
\label{RSF_subs}

\noindent
By a \sf{stratified space} (\sf{of dimension $k$}), we will mean a topological
space~$\ov\M$ together with a partition
$$\ov\M=\bigsqcup_{l=0}^{l=k}\M^{(l)}$$
such that $\M^{(l)}$ is a smooth manifold of (real) dimension $k\!-\!l$ and 
$$\ov\M^{(l)}-\M^{(l)}\subset\bigsqcup_{l'=l+1}^{l=k}\!\!\!\M^{(l')}\,.$$
If $U$ is an open subspace of a stratified space $\ov\M$ as above,
then 
$$U=\bigsqcup_{l=0}^{l=k}(\M^{(l)}\!\cap\!U)$$
is also a stratified space.
If $\ov\M_1$ and $\ov\M_2$ are stratified spaces, 
$\ov\M_1\!\times\!\ov\M_2$ is a stratified space with 
the strata given by unions of the products of the strata of 
$\ov\M_1$ and~$\ov\M_2$.
A continuous map $\pi\!:\ov\M_1\!\lra\!\ov\M_2$ between stratified spaces will
be called a \sf{stratified map} if the restriction of~$\pi$ to each
stratum of~$\ov\M_1$ is a smooth map to a stratum of~$\ov\M_2$.
A stratified map $\pi_V\!:V\!\lra\!\ov\M$ will be called a 
\sf{stratified vector bundle} if $\pi_V$ is a topological vector bundle
with fiber~$\C^k$ and the transition maps from open subsets of~$\ov\M$ to
$\GL_k\C$ are stratified.\\

\noindent
For the purposes of Definition~\ref{flatfam_dfn} below, we set
$$\pi_{\std}\!\equiv\!\pi_1: 
\fU_{\std}\equiv\big\{(t,u,v)\!\in\!\C^3\!:~ uv\!=\!t\big\}\lra\C$$
to be the projection to the first component.
This is a stratified map with respect to the stratifications
$$\C=\C^*\sqcup\{0\}, \qquad 
\fU_{\std}=\pi_{\std}^{-1}(\C^*)\sqcup \big(\pi_{\std}^{-1}(0)\!-\!0\big)
\sqcup\{0\}.$$
For each $t\!\in\!\C^*$, define
$$\rho_t\!:\Si_t\!\equiv\!\pi_{\std}^{-1}(t)\lra\R^+ \qquad\hbox{by}\quad 
\rho_t(t,u,v)=u^2+v^2\,.$$
If in addition $\ep\!\in\!\R^+$, let 
$$\Si_{t,\ep}=\big\{(t,u,v)\!\in\!\Si_t\!:\, |u|^2\!+\!|v|^2<\ep\big\}.$$
If $E\!\lra\!\Si_t$ is a normed vector bundle and $\eta\!\in\!L^p(\Si_t;E)$, let
$$\|\eta\|_{t,\ep}=\bigg(\int_{\Si_{t,\ep}}|\eta|^p\bigg)^{1/p}
+\bigg(\int_{\Si_{t,\ep}}\rho_t^{-\frac{p-2}{p}}|\eta|^2\bigg)^{1/2}\,.$$

\begin{dfn}\label{flatfam_dfn}
A stratified map $\pi\!:\fU\!\lra\!\ov\M$ is a 
\sf{flat stratified family of Riemann surfaces}~if
\begin{itemize}
\item each fiber $\Si_u\!\equiv\!\pi^{-1}(u)$ is a (possibly nodal) Riemann surface;
\item if $z_0\!\in\!\Si_{u_0}$ is a smooth point, 
there are neighborhoods $U_{z_0}$ of $u_0$ in $\ov\M$
and $\ti{U}_{z_0}$ of $z_0$ in~$\fU$ and a stratified isomorphism of fiber bundles 
$$\ti\phi_{z_0}\!:\ti{U}_{z_0}\lra U_{z_0}\!\times\!(\Si_{u_0}\!\cap\!\ti{U}_{z_0}) $$
over $U_{z_0}$ such that 
the restriction of~$\ti\phi_{z_0}$ to each fiber of~$\pi$ is holomorphic
and the restriction of~$\ti\phi_{z_0}$ to $\Si_{u_0}\!\cap\!\ti{U}_{z_0}$ is the identity;
\item if $z_0\!\in\!\Si_{u_0}$ is a node, there are neighborhoods $U_{z_0}$ of $u_0$ 
in $\ov\M$ and $\ti{U}_{z_0}$ of $z_0$ in~$\fU$, a stratified space~$U_{z_0}'$, and 
stratified embeddings
$$\phi_{z_0}\!: U_{z_0}\lra U_{z_0}'\times\C \qquad\hbox{and}\qquad
\ti\phi_{z_0}\!:\ti{U}_{z_0}\lra U_{z_0}'\times\!\fU_{\std} $$
such that the diagram 
$$\xymatrix{ \ti{U}_{z_0} \ar[d]^{\pi} \ar[r]^-{\ti\phi_{z_0}}& 
U_{z_0}'\!\times\!\fU_{\std} \ar[d]^{\id\times\pi_{\std}}\\
U_{z_0} \ar[r]^-{\phi_{z_0}}& U_{z_0}'\!\times\!\C }$$
commutes  and 
the restriction of~$\ti\phi_{z_0}$ to each fiber of~$\pi$ is holomorphic.
\end{itemize}
\end{dfn}

\begin{dfn}\label{flatfam_df2}
If $S$ is a finite set, a stratified map $\pi\!:\fU\!\lra\!\ov\M$ 
with stratified sections $z_j\!:\ov\M\!\lra\!\fU$, $j\!\in\!S$,
is a \sf{flat stratified family of $S$-marked Riemann surfaces}~if
\begin{itemize}
\item  $\pi\!:\fU\!\lra\!\ov\M$ 
is a flat stratified family of Riemann surfaces;
\item $z_j(u)\!\in\!\Si_u$ is a smooth point for every $u\!\in\!\ov\M$
and $j\!\in\!S$;
\item $z_{j_1}(u)\!\neq\!z_{j_2}(z)$ for every 
$u\!\in\!\ov\M$, $j_1,j_2\!\in\!S$ with $j_1\!\neq\!j_2$.
\end{itemize}
\end{dfn}

\begin{dfn}\label{flatfam_df3}
If $\pi\!:\fU\!\lra\!\ov\M$ is a flat stratified family of $S$-marked Riemann surfaces
and $Y$ is a smooth manifold,
a continuous map $F\!:\fU\!\lra\!Y$ is a \sf{flat family of $S$-marked maps} if
\begin{itemize}
\item for every $u\!\in\!\ov\M$, 
the restriction of $F$ to $\Si_u\!\equiv\!\pi^{-1}(u)$ is an $L^p_1$-map;
\item if $z_0\!\in\!\Si_{u_0}$ is a smooth point and $U_{z_0}$, $\ti{U}_{z_0}$,
and $\ti\phi_{z_0}$ are as in Definition~\ref{flatfam_dfn}, 
there exists a compact neighborhood $K_{z_0}(F)$ of $z_0$ in $\Si_{u_0}\!\cap\!\ti{U}_{z_0}$
such that $F\circ\ti\phi_{z_0}^{-1}|_{u\times K_{z_0}(F)}$ converges to
$F|_{K_{z_0}(F)}$ in the $L^p_1$-norm as $u\!\in\!U_{z_0}$ approaches~$u_0$;
\item if $z_0\!\in\!\Si_{u_0}$ is a node and $U_{z_0}$, $\ti{U}_{z_0}$,
$\phi_{z_0}$, and $\ti\phi_{z_0}$ are as in Definition~\ref{flatfam_dfn}, 
$$\lim_{\ep\lra0}\lim_{\underset{(u',t)\in\phi_{z_0}(U_{z_0})}{(u',t)\lra\phi_{z_0}(u)}} 
\big\|d(F\circ\ti\phi_{z_0}^{-1}|_{u'\times\Si_t})\big\|_{t,\ep}=0\,.$$
\end{itemize}
\end{dfn}

\noindent
In the case of interest to us, $\ov\M$ will be a family of 
$S$-marked stable maps to a smooth manifold~$Y$.
The fiber of $\fU\!\lra\!\ov\M$ over a point $u\!:\Si_u\!\lra\!Y$
will be the Riemann surface~$\Si_u$.

\subsection{Families of generalized CR-operators}
\label{CR_subs2}

\noindent
Let $D$ be a smooth CR-operator on a vector bundle $(N,\fI)$ over 
an almost complex manifold~$(Y,J)$.
Suppose $\fU\!\lra\!\ov\M$ is a flat stratified family of $S$-marked Riemann surfaces,
$F\!: \fU\!\lra\!Y$ is a flat family of maps, $S_0\!\subset\!S$, 
and $\vr$ is a function assigning to each $z_0\!\in\!S_0$
a real subbundle of $\ev_j^*N^*$.
For each $u\!\in\!\ov\M$ and $z_0\!\in\!S$, let $\vr_u(z_0)$ be the fiber
of $\vr(z_0)$ over~$u$.
Denote by $\ker_{\vr;u}^F(D)$ and $\ker_{\vr;u}^F(D^*)$
the kernels of the operators
\begin{equation*}\begin{split}
\big\{\big(F|_{\Si_u}\big)^*D\big\}_{\vr_u}\!:\,
& L^p_1\big(\Si_u;\{F|_{\Si_u}^*N\}(-\vr_u)\big)
\lra L^p\big(\Si_u;T^*\Si^{0,1}\!\otimes\!_{\C}F|_{\Si_u}^*N\big),\\
\big\{\big(F|_{\Si_u}\big)^*D\big\}_{\vr_u}^*\!:\,
& L^p_1\big(\Si_u;\K_{\Si_u}\!\otimes\!_{\C}\{F|_{\Si_u}^*N\}(\vr_u)\big)\\
&\qquad\qquad\qquad
\lra L^p\big(\Si_u;T^*\Si_u^{0,1}\!\otimes_{\C}\K_{\Si_u}\!\otimes\!_{\C}
\{F|_{\Si_u}^*N\}(\{z_j(u)\}_{j\in S_0})\big),
\end{split}\end{equation*}
respectively.\\

\noindent
We topologize the sets
$$\ker_{\vr}^F(D)\equiv\bigsqcup_{u\in\ov\M}\ker_{\vr;u}^F(D)
\qquad\hbox{and}\qquad
\ker_{\vr}^F(D^*)\equiv\bigsqcup_{u\in\ov\M}\ker_{\vr;u}^F(D^*)$$
by point-wise convergence on compact subsets of the complement of
the special (nodal and marked) points of the fiber.
In other words, suppose $u_r\!\in\!\ov\M$, $r\!\in\!\Z^+$, is a sequence converging
to $u_0\!\in\!\ov\M$ and $\xi_r\!\in\!\ker_{\vr;u_r}^F(D')$ for $r\!\in\!\bar\Z^+$,
where $D'\!=\!D,D^*$ and $\bar\Z^+\!=\!\{0\}\!\sqcup\!\Z$.
The sequence~$\{u_r\}$ \sf{converges to~$\xi_0$} if for every smooth point 
$z_0\!\in\!\Si_{u_0}$, with $z_0\!\neq\!z_j(u)$ for $j\!\in\!S$,
there exists a compact neighborhood~$K_{z_0}(F)$ as in Definition~\ref{flatfam_df3} 
such that $\xi_r\circ\ti\phi_{z_0}^{-1}|_{u_r \times K_{z_0}(F)}$ converges pointwise 
to~$\xi_0|_{K_{z_0}(F)}$.\\

\noindent
By Carleman Similarity Principle \cite[Theorem~2.2]{FHoS}, 
if the restriction of an element $\xi$ of 
$\ker_{\vr;u}^F(D')$ to an open subset of a component~$\Si_{u;i}$ of~$\Si_u$
vanishes, then the restriction of~$\xi$ to~$\Si_{u;i}$ is zero as well.
This implies that the above convergence topology on~$\ker_{\vr}^F(D)$ 
is the topology inherited from
the convergence topology on the bundle over $\ov\M$ with fibers
$L^p_1(\Si_u;u^*N)$ described in \cite[Section~3]{LiT}.\footnote{
While \cite[Section~3]{LiT} concerns only the case $N\!=\!TY$, 
it applies to any vector bundle $N\!\lra\!Y$.}
Furthermore, if the dimension of $\ker_{\vr;u}^F(D)$ is independent of~$u$,
then $\ker_{\vr}^F(D)\!\lra\!\ov\M$ is a vector bundle.
By \cite[Section~6]{RT}, the analogous statement holds 
for~$\ker_{\vr;u}^F(D^*)$.\footnote{While \cite[Section~6]{RT} concerns only 
the case $N\!=\!TY$ and $S_0\!=\!\eset$,
the argument applies to any vector bundle $N\!\lra\!Y$.
Furthermore, the twisting construction of \cite[Lemma~2.4.1]{Sh}
reduces the situation to the case $S_0\!=\!\eset$.
By \cite[Chapter~4]{Si}, which builds on \cite{SeSi}, 
there are Fredholm operators defining these vector spaces 
that form a continuous family  over $\ov\M$
and thus define a K-theory class; 
however, this statement is stronger than needed here.}
Lemma~\ref{serre_lmm} then implies that $\ker_{\vr}^F(D^*)\!\lra\!\ov\M$
is a vector bundle 
if the dimension of $\ker_{\vr;u}^F(D)$ is independent of~$u\!\in\!\ov\M$.
If in addition, the vector bundles $\ker_{\vr}^F(D)\!\lra\!\ov\M$ 
and $\vr(z_0)$, $z_0\!\in\!S$, are oriented 
(and $S$ is ordered if any of the bundles $\vr(z_0)$ is of odd rank),
then the vector bundle 
\BE{cokbnd_e} \ker_{\vr}^F(D^*)\lra \ov\M         \EE
has a canonical induced orientation,
since $\ker_{\vr;u}^F(D)$ and $(\ker_{\vr;u}^F(D^*))^*$ 
are the kernel and cokernel of an operator
obtained by a zeroth-order deformation from a first-order complex-linear Fredholm operator;
the determinant line of such an operator has a canonical orientation 
defined via a homotopy of  Fredholm operators
(see the proof of \cite[Theorem~3.1.5]{McS}).

\section{Proof of Theorem~\ref{main_thm}}
\label{GW_sec}

\noindent
The first claim of Theorem~\ref{main_thm} is immediate from the assumption
that $f_j^{-1}(Y)$ is a smooth oriented manifold.
Thus,
$$\big[\ov\M_{g,\f}(Y,\be_Y;J)\big]^{vir} =\bigg(
\prod_{j\in S}\big\{\ev_j\!\times\!(f_j\!\circ\!\pi_j)\!\big\}^*\big(\PD_{Y^2}(\De_Y)\big)
\bigg) \cap 
\bigg[\ov\M_{g,S}(Y,\be_Y;J)\times\prod_{j\in S}f_j^{-1}(Y)\bigg]^{vir}\,,$$
where $\De_Y\!\subset\!Y^2$ is the diagonal and 
$\pi_j\!:\prod_{j\in S}f_j^{-1}(Y)\lra f_j^{-1}(Y)$
is the projection onto the $j$-th component;
the identity~\e_ref{Ydim_e} now follows from~\e_ref{rkNf_e}.
Sections~\ref{CR_subs} and~\ref{CR_subs2} imply the second 
claim of Theorem~\ref{main_thm}.
Since the vector spaces 
\BE{kercok_e} \ker\big((D_{J;u}^{NY})^*\big)\approx \cok\big(D_{J;u}^{NY}\big)^* \EE 
have constant rank and are oriented via the isomorphism~\e_ref{mainthm_e}, 
they form natural oriented bundles over the uniformizing charts for 
$\ov\M_{g,\f}(Y,\be_Y;J)$ described in \cite[Section~3]{LiT}.
These bundles glue together to form an oriented vector orbi-bundle 
over $\ov\M_{g,\f}(Y,\be_Y;J)$.\footnote{Neither the topologies of the bundles
over the uniformizing charts nor the isomorphisms~\e_ref{kercok_e} depend 
on the Riemannian metrics over the uniformizing charts of \cite[Section~3]{LiT}.}
In the notation of Sections~\ref{CR_subs} and~\ref{CR_subs2}, 
this is also the bundle of the cokernels of the injective operators
$D_{J,\vr;\u}^{NY}\!\equiv\!(D_{J;u}^{NY})_{\vr}$, where 
\BE{udfn_e} [\u]\equiv \big([u],(w_j)_{j\in S}\big)\in \ov\M_{g,\f}(Y,\be_Y;J)\EE
and $\vr$ is the function
assigning to each element $j\!\in\!S$ the subbundle 
$\Ann(\ev_j^*(\Im\, d^{NY}f_j),\R)$ of~$\ev_j^*NY^*$.
The identity~\e_ref{cokrk_e} is immediate from~\e_ref{mainthm_e} and the Index Theorem.
The first part of the third claim follows immediately from
Proposition~\ref{horreg_prp} below in light of assumption~(b)
in Theorem~\ref{main_thm}.\\

\noindent
We note that the second part of the third claim of Theorem~\ref{main_thm}
is consistent with the divisor relation for GW-invariants \cite[(3.4)]{RT2}
in the following sense.
Let
$$\pi_0\!: \ov\M_{g,\{0\}\sqcup S}(Y,\be_Y;J)\lra \ov\M_{g,S}(Y,\be_Y;J)$$
be the forgetful map dropping the 0-th marked point and 
$f_0\!:M_0\!\lra\!X$ a cobordism representative for some $\ka_0\!\in\!H_{2n-2}(X;\Z)$
so that $f_0$ is transverse to~$Y$;
the last assumption implies that $N_{f_0(w_0)}^{f_0}Y\!=\!\{0\}$ for 
all $w_0\!\in\!f_0^{-1}(Y)$.
With 
\begin{equation*}\begin{split}
\ov\M_{g,f_0\sqcup\f}(Y,\be_Y;J)
\equiv\big\{\big([u],w_0,(w_j)_{j\in S}\big)\in
\ov\M_{g,\{0\}\sqcup S}(Y,\be_Y;J)\!\times\!M_0\!\times\!\prod_{j\in S}\!M_j\!:
\qquad\qquad\qquad&\\
\ev_j([u])\!=\!f_j(w_j)~\forall\,j\!\in\!\{0\}\!\sqcup\!S\big\}\,,&
\end{split}\end{equation*}
let 
$$\ti\pi_0\!:\ov\M_{g,f_0\sqcup\f}(Y,\be_Y;J)\lra \ov\M_{g,\f}(Y,\be_Y;J)$$
be the map induced by $\pi_0$.
If $[u]\!\in\!\ov\M_{g,\{0\}\sqcup S}(Y,\be_Y;J)$
and $\pi_0$ contracts component $\Si_{u;i_0}$ of $\Si_u$, then 
$\Si_{u;i_0}$ is $\P^1$,  contains precisely two nodes, 
say $0$ and~$\i$, along with the $0$-th marked point and no other marked points,
and $u|_{\Si_{u;i_0}}$ is constant.
Therefore, if $[u']\!=\!\pi_0(u)$ and $\chi_u$ is the set of components of $\Si_u$, 
then the homomorphisms
\begin{alignat}{2}
\label{kerisom_e}
&\ker(D_{J;u}^{NY})\lra  \ker(D_{J;u'}^{NY}), &\qquad 
&(\xi_i)_{i\in\chi_u}\lra (\xi_i)_{i\in\chi_u-i_0}\,,\\
\label{cokisom_e}
&\ker\big((D_{J;u}^{NY})^*\big)\lra  \ker\big((D_{J;u'}^{NY})^*\big), &\qquad 
&(\eta_i)_{i\in\chi_u}\lra (\eta_i)_{i\in\chi_u-i_0}\,,
\end{alignat}
are well-defined and are in fact isomorphisms.
Since \e_ref{kerisom_e} is an isomorphism, $f_0\!\sqcup\!\f$ satisfies
the assumptions of Theorem~\ref{main_thm} if and only if $\f$ does.
Since the total spaces of the cokernel bundles are topologized using 
convergence of elements of $\ker(D_{J;u}^{NY})^*$ 
on compact subsets of smooth points,
\e_ref{cokisom_e} induces an isomorphism of orbi-bundles
\BE{cokisom_e2} \cok(D_J^{NY})\lra \ti\pi_0^*\cok(D_J^{NY})\EE 
over $\ov\M_{g,f_0\sqcup\f}(Y,\be_Y;J)$;
it extends over a neighborhood of $\ov\M_{g,f_0\sqcup\f}(Y,\be_Y;J)$
in the space of $L^p_1$-maps via the construction described at the 
end of Section~\ref{config_subs}.
Thus, by the standard divisor relation, 
\begin{equation*}\begin{split}
&\bigg\lan e\big(\cok(D_J^{NY})\big)\prod_{j\in S}\psi_j^{a_j},
\big[\ov\M_{g,f_0\sqcup\f}(Y,\be_Y;J)\big]^{vir}\bigg\ran\\
&\hspace{1in}
=\blr{\PD_Y\ka_0,\be}\cdot
\bigg\lan e\big(\cok(D_J^{NY})\big)\prod_{j\in S}\psi_j^{a_j},
\big[\ov\M_{g,\f}(Y,\be_Y;J)\big]^{vir}\bigg\ran.
\end{split}\end{equation*}
In particular, it is sufficient to verify \e_ref{mainthm_e2}
under the assumption that $2g\!+\!|S|\!\ge\!3$; 
this slightly simplifies the presentation.\\

\noindent
For the remainder of the paper, we assume that $2g\!+\!|S|\!\ge\!3$.
Section~\ref{config_subs} sets up notation for the configuration spaces
that play a central role in \cite{FuO} and~\cite{LiT}.
The main geometric observation used in the proof of Theorem~\ref{main_thm}
is Proposition~\ref{horreg_prp}, stated and proved in Section~\ref{subman_subs}.
Our approach to~\e_ref{mainthm_e2} is illustrated in Section~\ref{semipos_subs},
where \e_ref{mainthm_e2} is verified in some cases, including
the case of Theorem~\ref{FanoGV_thm}.
The general case is the subject of Section~\ref{mainpf_subs}.

\subsection{Configuration spaces}
\label{config_subs}

\noindent
Let $X$ be a compact manifold, $\be\!\in\!H_2(X;\Z)$, $g$ a non-negative integer,
and $S$ a finite set. 
We denote by $\X_{g,S}(X,\be)$  the space of equivalence classes
of stable $L^p_1$-maps $u\!:\Si_u\!\lra\!X$ from genus~$g$  Riemann surfaces
with $S$-marked points, which may have simple nodes, to~$X$ of degree~$\be$, i.e. 
$$u_*[\Si_u]=\be\in H_2(X;\Z).$$
Let $\X_{g,S}^0(X,\be)$ be the subset of $\X_{g,S}(X,\be)$
consisting of the stable maps with smooth domains. 
The space $\X_{g,S}(X,\be)$ is topologized in \cite[Section~3]{LiT}
using $L^p_1$-convergence on compact subsets 
of smooth points of the domain and certain convergence requirements near the nodes.
The space $\X_{g,S}(X,\be)$ is stratified by subspaces $\X_{\T}(X)$ of stable maps 
from domains of the same geometric type and with
the same degree distribution between the components of the domain.
Each stratum is the quotient of a smooth Banach manifold 
$\ti\X_{\T}(X)$ by a finite-dimensional Lie group~$G_{\T}$;
the restriction of the $G_{\T}$-action to any finite-dimensional submanifold
of $\ti\X_{\T}(X)$ consisting of smooth maps and preserved by $G_{\T}$ is smooth.
The closure of the main stratum, $\X_{g,S}^0(X,\be)$, is $\X_{g,S}(X,\be)$.
If  $f_j\!:M_j\!\lra\!X$ for $j\!\in\!S$ are smooth maps, let
$$\X_{g,\f}(X,\be)=\big\{\big([u],(w_j)_{j\in S}\big)
\in \X_{g,S}(X,\be)\times\prod_{j\in S}M_j\!:~
u(z_j(u))\!=\!f_j(w_j)~\forall\,j\!\in\!S\big\}.$$\\

\noindent
If $J$ is an almost complex structure on $X$, let 
$$\Ga_{g,S}^{0,1}(X,\be;J)\!\lra\!\X_{g,S}(X,\be)$$
be the family of $(TX,J)$-valued $(0,1)$ $L^p$-forms. 
In other words, the fiber of $\Ga_{g,S}^{0,1}(X,\be;J)$ over a point
$[u]$ in $\X_{g,S}(X,\be)$ is the space
$$\Ga_{g,S}^{0,1}(X,\be;J)\big|_{[u]}=\Ga^{0,1}(X,u;J)\big/\hbox{Aut}(u),
\quad\hbox{where}\quad
\Ga^{0,1}(X,u;J)=L^p\big(\Si_u;T^*\Si_u^{0,1}\!\otimes\!_{\C}u^*TX\big).$$
The total space of this family is topologized in \cite[Section~3]{LiT}
using $L^p$-convergence on compact subsets of smooth points of the domain
and certain convergence requirements near the nodes.
The restriction of $\Ga_{g,S}^{0,1}(X,\be;J)$ to each stratum $\X_{\T}(X)$ is 
the quotient of a smooth Banach vector bundle $\ti\Ga_{\T}^{0,1}(X;J)$
over $\ti\X_{\T}(X)$ by~$G_{\T}$.
The smooth sections of the bundles  
$\ti\Ga_{\T}^{0,1}(X;J)\lra \ti\X_{\T}(X)$
given~by
$$\dbar_J\big([\Si_u,\fJ_u;u]\big) = \dbar_{J,\fJ_u}u
= \frac{1}{2}\big(du+J\!\circ\!du\!\circ\!\fJ_u\big)$$
induce sections of $\Ga_{g,S}^{0,1}(X,\be;J)$ over $\X_{\T}(X)$,
which define a continuous section $\bar\partial_J$ of the family 
$$\Ga_{g,S}^{0,1}(X,\be;J) \lra \X_{g,S}(X,\be).$$ 
The zero set of this section
is the moduli space $\ov\M_{g,S}(X,\be;J)$ 
of equivalence classes of stable $J$-holomorphic degree~$\be$ maps 
from genus-$g$ curves with $S$-marked points into~$X$.
The section $\dbar_J$ over $\ti\X_{\T}(X)$ is Fredholm, 
i.e.~its linearization has finite-dimensional kernel and cokernel
at every point of the zero set.
The index of the linearization $D_{J;u}$ of~$\dbar_J$ at $u\!\in\!\ti\X_{\T}(X)$
such that 
$$[u]\in\M_{g,S}(X,\be;J)\equiv \ov\M_{g,S}(X,\be;J)\cap \X_{g,S}^0(X,\be)$$ 
is the expected dimension $\dim_{g,S}(X,\be)$ of the moduli space 
$\ov\M_{g,S}(X,\be;J)$. \\

\noindent
If $f_j\!:M_j\!\lra\!X$ for $j\!\in\!S$ are smooth maps, 
$Y\!\subset\!X$ is a submanifold, $\be_Y\!\in\!H_2(X;\Z)$ is such that 
$\io_{Y*}\be_Y\!=\!\be$,
and $\T$ is any combinatorial type of maps to $X$ or $Y$
of degree $\be$ or $\be_Y$, respectively, let
\begin{equation*}\begin{split}
\X_{g,\f}(X,\be)&=\big\{\big([u],(w_j)_{j\in S}\big)\in 
\X_{g,S}(X,\be)\!\times\!\prod_{j\in S}\!M_j\!:\,
\ev_j([u])\!=\!f_j(w_j)\,\forall\,j\!\in\!S\big\},\\
\X_{g,\f}(Y,\be_Y)&=\X_{g,\f}(X,\be)
\cap \bigg(\X_{g,S}(Y,\be_Y)\!\times\!\prod_{j\in S}\!M_j\bigg),\\
\X_{\T,\f}(X)&=\X_{g,\f}(X,\be)
\cap \bigg(\X_{\T}(X)\!\times\!\prod_{j\in S}\!M_j\bigg),\\
\X_{\T,\f}(Y)&=\X_{g,\f}(Y,\be_Y)\cap 
\bigg(\X_{\T}(Y)\!\times\!\prod_{j\in S}\!M_j\bigg).
\end{split}\end{equation*}
With $\pi\!:\X_{g,\f}(X,\be)\!\lra\!\X_{g,S}(X,\be)$ denoting 
the projection map, let
\begin{equation*}\begin{split}
\Ga^{0,1}_{g,\f}(X,\be;J)&=\pi^*\Ga^{0,1}_{g,S}(X,\be;J)\lra \X_{g,\f}(X,\be);\\
\Ga^{0,1}_{g,\f}(Y,\be_Y;J)&=\pi^*\Ga^{0,1}_{g,S}(Y,\be_Y;J)\lra \X_{g,\f}(Y,\be_Y).
\end{split}\end{equation*}
With $a_j$, $j\!\in\!S$, as in Theorem~\ref{main_thm}, let
$$\bL_{\a,\f}\equiv\bigoplus_{j\in S}a_j\pi^*L_j^*\lra\X_{g,\f}(X,\be),$$
where $L_j\!\lra\!\X_{g,S}(X,\be)$ is the tautological line bundle
for the $j$-th marked point.\\

\noindent
If $J$ is an almost complex structure on $X$ preserving $Y$, 
let $g_J$ be a $J$-invariant metric on $X$, $\na^J$ the $J$-linear connection
of $g_J$ induced by the Levi-Civita connection of $g_J$, 
$TY^{\v}\!\subset\!TX|_Y$ the $g_J$-orthogonal complement of $TY$,
and $\pi^{\h}\!:TX|_Y\!\lra\!TY$ the orthogonal projection map.
Define
\begin{gather*}
\ti\na^J\!: \Ga(Y;TX)\lra \Ga(Y;T^*Y\!\otimes_{\R}\!TX)  \qquad\hbox{by}\\
\ti\na^J_v(\xi^{\h}+\xi^{\v})=
\pi^{\h}\big(\na_v^J\xi^{\h}\big)+\na_v^J\xi^{\v}
\qquad\forall~v\!\in\!TY,~\xi^{\h}\!\in\!\Ga(Y;TY),~\xi^{\v}\!\in\!\Ga(Y;TY^{\v}).
\end{gather*}
This connection in $TX|_Y$ gives rise to a $\C$-linear connection $\na^{\perp}$ on $NY$
and thus to a $\dbar$-operator $\dbar^{\perp}$ on~$NY$.
Define
$$D^{NY}\!:\Ga(Y;NY)\lra \Ga(Y;T^*Y^{0,1}\!\otimes_{\C}\!NY)  
\qquad\hbox{by}\qquad
D^{NY}\xi=\dbar^{\perp}\xi+A_X^{\perp}(\cdot,\xi),$$
where $A_X^{\perp}$ is the composition of the Nijenhuis tensor of $J$ on $X$ 
with the projection to~$NY$.
If $[u]\!\in\!\X_{g,S}(Y,\be_Y)$, let 
$$D^{NY}_{J;u}\!: L^p_1(\Si_u;u^*NY)\lra 
L^p\big(\Si_u;T^*\Si_u^{0,1}\!\otimes_{\C}\!u^*NY\big)$$
be the pull-back of $D^{NY}$ by $u$ with respect to the connection $\na^{\perp}$
as in Section~\ref{CR_subs}.
If $[u]$ is an element of $\ov\M_{g,S}(Y,\be_Y;J)$, 
this definition agrees with the one in Section~\ref{AbsGW_subs}.
Thus, under the assumptions of Theorem~\ref{main_thm}, 
the dimension of $\cok(D^{NY}_{J;u})$ is fixed on a neighborhood of 
$\ov\M_{g,\f}(Y,\be_Y;J)$ in $\X_{g,\f}(Y,\be_Y)$.
By Section~\ref{CR_subs2}, the vector spaces $\cok(D^{NY}_{J;u})$  form a vector 
orbi-bundle over such a neighborhood.

\subsection{Symplectic submanifolds and pseudo-holomorphic maps}
\label{subman_subs}

\begin{dfn}\label{Jtub_dfn}
If $(X,J)$ is an almost complex manifold and $Y\!\subset\!X$ is 
an almost complex submanifold, a tuple
$(\pi_Y\!:U_Y\!\lra\!Y,TU_Y^h)$ 
is a \sf{$J$-regularized tubular neighborhood of~$Y$ in~$X$}~if 
\begin{itemize}
\item $U_Y$ is a tubular neighborhood of~$Y$ in~$X$;
\item $\pi_Y\!:U_Y\!\lra\!Y$ is a vector bundle such that $\pi_Y|_Y\!=\!\id_Y$ and
$\ker d_y\pi_Y$ is a complex subspace of $(T_yX,J)$ 
for every $y\!\in\!Y$;
\item $TU_Y^h\!\lra\!U_Y$ is a complex subbundle of $(TU_Y,J)$ such that 
$d_x\pi_Y\!:TU_Y^h\!\lra\!T_{\pi_Y(x)}Y$ is an isomorphism 
of real vector spaces for every $x\!\in\!U_Y$
and is the identity for every $x\!\in\!Y$.\\
\end{itemize}
\end{dfn}

\noindent
Every embedded almost complex submanifold $Y$ of an almost complex
manifold $(X,J)$ admits a $J$-regularized tubular neighborhood.
Let $g$ be a $J$-invariant Riemannian metric on $X$ and $\exp^g\!:TX\!\lra\!X$ 
the exponential map with respect to the Levi-Civita connection of the metric~$g$.
Identifying $NY$ with the $g$-orthogonal complement of $TY$ in $TX|_Y$, 
we obtain a smooth map
$$\exp^Y\!: NY\lra X$$
by restricting $\exp^g$.
Since $Y$ is an embedded submanifold of $X$, there exist 
tubular neighborhoods $U_Y'$ and $U_Y$ of $Y$ in $NY$ 
and in $Y$, respectively, such that the map
$$\exp\!\equiv\!\exp^Y\big|_{U_Y'}\!: U_Y'\lra U_Y$$
is a diffeomorphism. Furthermore, $\exp|_Y\!=\!\id_Y$ and 
$d_y\exp\!:T_yNY\lra T_yX$ is $\C$-linear for every $y\!\in\!Y$.
Thus, 
$$\pi_Y=\pi_{NY}\!\circ\!\exp|_{U_Y'}^{\,-1}\!: U_Y\lra Y,$$
where $\pi_{NY}\!: NY\!\lra\!Y$ is the bundle projection map,
satisfies the middle condition in Definition~\ref{Jtub_dfn}.
Furthermore, if $(\ker d\pi_Y)^{\perp}$ is the $g$-orthogonal complement
of $\ker d\pi_Y$ in~$TU_Y$,
$$d_x\pi_Y\!:(\ker d_x\pi_Y)^{\perp}\lra T_{\pi_Y(x)}Y$$
is an isomorphism and induces a complex structure $J_Y$
in the vector bundle $(\ker d\pi_Y)^{\perp}\!\lra\!U_Y$
(which may differ from~$J$).
Let
$$T_xU_Y^h=\big\{v\!-\!JJ_Yv\!: v\!\in\!(\ker d_x\pi_Y)^{\perp}\big\}.$$
Note that $T_xU_Y^h$ is a complex linear subspace of $(T_xU_Y,J_x)$
for each $x\!\in\!U_Y$.
Since $(\ker d_y\pi_Y)^{\perp}\!=\!T_yY$ and
$J_Y|_y\!=\!J|_{T_yY}$ for every $y\!\in\!Y$, 
$$d_y\pi_Y=\id\!: T_yU_Y^h\lra T_{\pi_Y(y)}Y$$
for every $y\!\in\!Y$. 
Thus,
$$d_x\pi_Y\!: T_xU_Y^h\lra T_{\pi_Y(x)}Y$$
is an isomorphism for every $x\!\in\!U_Y$ if $U_Y$ is sufficiently small.
We conclude that $TU_Y^h$ satisfies the final condition in 
Definition~\ref{Jtub_dfn}.

\begin{prp}\label{horreg_prp}
Suppose $(X,\om)$ is a compact symplectic manifold, $g\!\in\!\bar\Z^+$,
$S$ is a finite set, $\be\!\in\!H_2(X;\Z)$,  
and $f_j\!:M_j\!\lra\!X$ is a smooth map for each $j\!\in\!S$.
Let $J$ be an $\om$-tame almost complex structure on~$X$,
$Y$ a compact almost complex submanifold of~$(X,J)$, and 
$(\pi_Y\!:U_Y\!\lra\!Y,TU_Y^h)$ 
a $J$-regularized tubular neighborhood of~$Y$ in~$X$.
If $([u_r],(w_{r,j})_{j\in S})\in\X_{g,\f}(X,\be)$ is a sequence such that
\begin{gather}\label{dbarhor_e}
u_r(\Si_{u_r})\not\subset Y, \qquad
\dbar_Ju_r\big|_{u_r^{-1}(U_Y)}\in 
L^p\big(u_r^{-1}(U_Y);T^*(u_r^{-1}(U_Y))^{0,1}\!\otimes_{\C}\!u_r^*TU_Y^h\big),\\
\lim_{r\lra\i}\big([u_r],(w_{r,j})_{j\in S})=\big([u],(w_j)_{j\in S})
\in\ov\M_{g,\f}(Y,\be_Y;J)\subset\X_{g,\f}(X,\be)
\notag
\end{gather}
for some $\be_Y\!\in\!H_2(Y;\Z)$,  then
$$\exists~~ \xi\!\in\!\ker\,D_{J;u}^{NY}, ~~ v_j\!\in\!T_{w_j}M_j~
\forall\,j\!\in\!S
\qquad\hbox{s.t.}\qquad 
\xi\neq0, \qquad \xi\big(z_j(u)\big)=d_{w_j}f_j(v_j) \quad\forall\,j\!\in\!S.$$\\
\end{prp}

\noindent
The rest of this subsection is dedicated to the proof of this proposition
by adopting a now-standard rescaling argument.
It is sufficient to consider the case $X\!=\!NY$ as smooth manifolds
and $\pi_Y\!:NY\!\lra\!Y$ is the bundle projection map.
After passing to a subsequence, it can be assumed that 
the topological types of the domains~$\Si_{u_r}$ of~$u_r$ are the same
(but not necessarily the same as the topological type of~$\Si_u$). 
The desired vector field~$\xi$ and tangent vectors~$v_j$ will be constructed
by re-scaling~$u_r$ in the normal direction to~$Y$ and then taking the limit.\\

\noindent
For each $j\!\in\!S$, let $N_jY\!\subset\!T_{w_j}M$ be a complement 
of $T_{w_j}(f_j^{-1}(Y))$ and 
$$\exp_j\!: T_{w_j}M_j\lra M_j$$
a diffeomorphism onto a neighborhood of~$w_j$ in~$M_j$ such that 
$$\exp_j(0)=w_j, \qquad d_0\exp_j=\Id, \qquad 
\exp_j(v)\in f_j^{-1}(Y)~~\forall\,v\!\in\!T_{w_j}(f_j^{-1}(Y)).$$
For each $r\!\in\!\Z^+$, define
$$v_{r,j}^h\oplus v_{r,j}^{\perp}\in 
T_{w_j}(f_j^{-1}(Y))\oplus N_jY=T_{w_j}M_j
\qquad\hbox{by}\quad \exp_j\big(v_{r,j}^h\!+\!v_{r,j}^{\perp}\big)=w_{r,j}\,.$$
Choose metrics on~$NY$ and $N_jY$, $j\!\in\!S$.
By our assumptions, 
$$\ep_r\equiv\sup_{z\in\Si_{u_r}}\!\!\!\big|u_r(z)\big|\in\R^+,
\quad \lim_{r\lra\i}\!\!\ep_r=0,   \quad
\lim_{r\lra\i}\!v_{r,j}^h=0~~\forall\,j\!\in\!S, \quad
\big|v_{r,j}^{\perp}\big|\le C\ep_r~~\forall\,r\!\in\!\Z^+,\,j\!\in\!S,$$
for some $C\!\in\!\R^+$ independent of $r$ and $j$.
By the last condition, for each $j\!\in\!S$ (a subsequence of) the sequence
$$\ti{v}_{r,j}^{\perp}=\ep_r^{-1}v_{r,j}^{\perp}, \quad r\!\in\!\Z^+,$$
converges to some $v_j\!\in\!N_jY\!\subset\!T_{w_j}M_j$.\\

\noindent
For each $r\!\in\!\Z^+$, we define
\begin{alignat*}{2}
&m_r\!: NY\lra NY &\qquad &\hbox{by}\quad 
m_r(x)=\ep_r\cdot x;\\
&J_r\in\Ga\big(NY;\Hom(T(NY),T(NY))\big)&\qquad &\hbox{by}\quad
J_r|_x= \big\{d_x m_r\big\}^{-1}\circ J_{\ep_rx}\circ d_x m_r;\\
&\ti{u}_r\!: \Si_{u_r}\lra NY&\qquad &\hbox{by}\quad 
\ti{u}_r(z)=\ep_r^{-1}\cdot u_r(z);\\
&\eta_r\in L^p(\Si_{u_r};T^*\Si_{u_r}^{0,1}\!\otimes\!_{\C}\ti{u}_r^*T(NY))&\qquad 
&\hbox{by}\quad
\eta_r= \big\{d_{u_r(\cdot)}m_r\big\}^{-1}\circ\dbar_Ju_r.
\end{alignat*}
If in addition $j\!\in\!S$, define $\ti{f}_{r,j}\!:T_{w_j}M_j\!\lra\!NY$ by 
$$\ti{f}_{r,j}\big(v^h+v^{\perp}\big)=
\ep_r^{-1}\cdot f_j\big(\exp_j(v^h+\ep_rv^{\perp})\big)
\quad\forall\,v^h\!\in\!T_{w_j}(f_j^{-1}(Y)),\,v^{\perp}\!\in\!N_jY.$$
Then, for all $r\!\in\!\Z^+$,
\BE{rescalecond_e}\dbar_{J_r}\ti{u}_r=\eta_r, \quad
\sup_{z\in\Si_{\ti{u}_r}}\big|\ti{u}_r(z)\big|=1, \quad
\ti{u}_r(z_j(u_r))=\ti{f}_{r,j}\big(v_{r,j}^h+\ti{v}_{r,j}^{\perp}\big)
~~\forall\,j\!\in\!S.\EE
By the following paragraph, the sequence of almost complex structures 
$J_r$ $C^{\infty}$-converges on compact subsets of~$NY$ to an almost complex 
structure~$\ti{J}$ such that $\ti{J}|_{TY}\!=\!J|_{TY}$ 
and 
$$\dbar_{\ti{J}}\xi=0 ~~\Llra~~
D_{J;u}^{NY}\xi=0 \qquad\forall~\xi\in\Ga\big(\Si_u;u^*NY).$$
Furthermore, the sequence $\eta_r$ converges to~$0$.
Thus, by~\e_ref{rescalecond_e}, $\ti{u}_r$ converges to some 
\begin{gather*}
[\ti{u}]\in \ov\M_{g,S}(NY,\be;\ti{J})\subset \X_{g,S}(NY,\be)
\qquad\hbox{s.t.}\\
\ti{u}(\Si_{\ti{u}})\not\subset Y, ~~~
\ti{u}(x_j(\ti{u}))=d_{w_j}f_j(v_j)\in N_{f_j(w_j)}Y~~\forall\,j\!\in\!S.
\end{gather*}
Since we must have $\pi_Y\!\circ\!\ti{u}=\!u$, 
$\ti{u}$ corresponds to a section $\xi$ of $u^*NY\!\lra\!\Si_u$ as needed.\\ 

\noindent
It remains to prove the two local claims made above.
It is sufficient to assume that
$$\pi_Y\!=\!\pi_1\!: NY=Y\!\times\!\C^k\lra Y$$ 
as vector bundles over $Y$, and there exists
\begin{gather}
\al\in\Ga(Y\!\times\!\C^k;\Hom_{\R}(\pi_1^*TY,\pi_2^*T\C^k)\big)
\qquad\st \notag\\
\al|_{Y\times0}=0, \qquad 
\label{Jregulr_e2}
T_{(y,w)}U_Y^h=\big\{\big(y',\al_{(y,w)}(y')\big)\!: y'\!\in\!T_yY\big\}
\qquad \forall~(y,w)\in Y\!\times\!\C^k.
\end{gather}
Thus, by assumption on $u_r$,
$$\dbar_Ju_r=(\nu^{\h},\al_u\nu^h)
\qquad\hbox{for some}\quad \nu^h\in 
L^p(\Si_{u_r};T^*\Si_{u_r}\!\otimes_{\R}\!u_r^{\h*}TY\big),$$
where $u_r^{\h}=\pi_1\circ u_r$.
Let 
$$J=\left(\begin{array}{cc} J^{\h\h}& J^{\h\v}\\ J^{\v\h}& J^{\v\v}\end{array}\right)\!:
TU_Y\!=\!\pi_1^*TY\!\oplus\!\pi_2^*T\C^k\lra \pi_1^*TY\!\oplus\!\pi_2^*T\C^k $$
be the almost complex structure.
By Definition~\ref{Jtub_dfn},  $J^{\h\v}|_{Y\times0}\!=\!0$ and $J^{\v\h}|_{Y\times0}\!=\!0$;
we can also assume that $J^{\v\v}|_{Y\times0}\!=\!\fI$ is the standard complex 
structure on~$\C^k$.
If $\vec\na$ is the gradient with respect to the standard coordinates on~$\C^k$,
it follows that 
\begin{gather*}
\al_{(y,w)}=\ti\al_{(y,w)}w, \qquad J^{\v\h}_{(y,w)}=\ti{J}^{\v\h}_{(y,w)}w, \qquad 
J^{\v\v}_{(y,w)}=\fI+\ti{J}^{\v\v}_{(y,w)}w,  \qquad\hbox{where}\\
\ti\al_{(y,w)}=\int_0^1\vec\na \al_{(y,tw)}\,dt, \qquad
\ti{J}^{\v\h}_{(y,w)}=\int_0^1\vec\na J^{\v\h}_{(y,tw)}\,dt, \qquad
\ti{J}^{\v\v}_{(y,w)}=\int_0^1\vec\na J^{\v\v}_{(y,tw)}\,dt.\notag
\end{gather*}
This gives
\begin{equation*}\begin{split}
\eta_r&=\left(\begin{array}{c}\nu^{\h}\\ \ep_r^{-1}\{\ti\al_{u_r}u_r\}\nu_r^{\h}
\end{array}\right)\lra0\,,\\
J_r|_{(y,w)}&=
\left(\begin{array}{cc} J_{(y,\ep_rw)}^{\h\h}& \ep_r J_{(y,\ep_rw)}^{\h\v}\\ 
\ep_r^{-1} J_{(y,\ep_rw)}^{\v\h}& J_{(y,\ep_rw)}^{\v\v}\end{array}\right)
\lra \left(\begin{array}{cc} J_{T_yY}& 0\\ 
\ti{J}_{(y,0)}^{\v\h}w& \fI\end{array}\right)\equiv\ti{J}_{(y,w)}\,,\\
D_{J;u}\left(\begin{array}{c}\xi^{\h}\\ \xi^{\v}\end{array}\right)
&=\left(\begin{array}{c}\dbar\xi^{\h}\\ 
\dbar\xi^{\v}+\frac{1}{2}\{\ti{J}_{(y,0)}^{\v\h}\xi^{\v}\}du\circ\fJ\end{array}\right);
\end{split}\end{equation*}
the last identity is a special case of \cite[(3.1.4)]{McS}.
This concludes the proof of Proposition~\ref{horreg_prp}.

\subsection{Geometric motivation for \e_ref{mainthm_e2}}
\label{semipos_subs}

\noindent
In this section we give a rough argument for \e_ref{mainthm_e2}
before translating it into the virtual setting of \cite{FuO} and 
\cite{LiT} in Section~\ref{mainpf_subs}.
As explained at the end of this section, this argument suffices in some cases.
We continue with the notation of Theorem~\ref{main_thm} 
and Section~\ref{config_subs}.
For the remainder of the paper, we assume that \e_ref{dimcond_e} holds;
otherwise, the left-hand side of \e_ref{mainthm_e} vanishes by definition,
while the right-hand side vanishes by \e_ref{Ydim_e} and \e_ref{cokrk_e}.
Our assumption implies that 
\BE{dimeq_e} 
\dim_{g,\f}(Y,\be_Y)\equiv \big[\ov\M_{g,\f}(Y,\be_Y;J)\big]^{vir}
=2\sum_ja_j+\rk_{\R}\cok(D_J^{NY}).\EE
We also assume that $a_j\!\ge\!0$ for every $j\!\in\!S$.\\

\noindent
If $\nu$ is a sufficiently small multi-section of $\Ga^{0,1}_{g,\f}(X,\be;J)$
over $\X_{g,\f}(X,\be)$, the space
$$\ov\M_{g,\f}(X,\be;J,\nu)
=\{\dbar_J\!+\!\nu\}^{-1}(0)\subset \X_{g,\f}(X,\be)$$
is compact, because $\ov\M_{g,\f}(X,\be;J)$ is.
If in addition $\nu$ is smooth and generic in the appropriate sense, 
$\ov\M_{g,\f}(X,\be;J,\nu)$ is stratified by smooth branched orbifolds
of even dimensions.
If $\vph$ is a multi-section of the orbi-bundle $\bL_{\a,\f}\lra\X_{g,\f}(X,\be)$,
let
\BE{vhp0dfn_e}
\ov\M_{g,\f}^{\vph}(X,\be;J,\nu)=\ov\M_{g,\f}(X,\be;J,\nu)\cap \vph^{-1}(0).\EE
If $\nu$ is sufficiently small and generic and $\vph$ is generic,
the left-hand side of~\e_ref{mainthm_e2} is the number of elements of
$\ov\M_{g,\f}^{\vph}(X,\be;J,\nu)$ counted with appropriate multiplicities 
that lie in a small neighborhood of 
$$\ov\M_{g,\f}^{\vph}(Y,\be_Y;J)\equiv 
\ov\M_{g,\f}(Y,\be_Y;J)\cap\vph^{-1}(0)$$
in $\X_{g,\f}(X,\be)$.\\

\noindent
In order to verify \e_ref{mainthm_e2}, fix 
a $J$-regularized tubular neighborhood $(\pi_Y\!:U_Y\!\lra\!Y,TU_Y^h)$.
We will take $\nu\!=\!\nu_Y\!+\!\nu_X$ so that 
\begin{itemize}
\item for every $\u\!=\!([u],(w_j)_{j\in S})\!\in\!\X_{g,\f}(X,\be)$
with $[u]\!\in\!\X_{g,S}(U_Y,\be_Y)$,
$$\nu_Y(\u)\in L^p(\Si_u;T^*\Si_u^{0,1}\!\otimes_{\C}\!TU_Y^h);$$
\item $\nu_Y|_{\X_{g,\f}(Y,\be_Y)}$ is generic,
so that $\ov\M_{g,\f}(Y,\be_Y;J,\nu_Y)$ is stratified by smooth 
branched manifolds of the expected dimensions and the dimension of the main stratum
$$\M_{g,\f}(Y,\be_Y;J,\nu_Y) \equiv \ov\M_{g,\f}(Y,\be_Y;J,\nu_Y)
\cap \bigg(\X_{g,S}^0(Y,\be_Y)\times\prod_{j\in S}\!M_j\bigg)$$
is $\dim_{g,\f}(Y,\be_Y)$;
\item $\nu_X$ is generic and small relative to $\nu_Y$. 
\end{itemize} 
Using $\pi_Y$, $d\pi_Y|_{TU_Y^h}^{\,-1}$, and a bump function around $Y$ 
with support in $U_Y$, any section of 
$$\pi^*\Ga^{0,1}_{g,S}(Y,\be_Y;J)\lra
\X_{g,S}(Y,\be_Y)\times\prod_{j\in S}\!M_j$$
can be extended
to a section of  $\Ga^{0,1}_{g,\f}(X,\be;J)$ over $\X_{g,\f}(X,\be)$
satisfying the middle condition above.
In light of Proposition~\ref{horreg_prp}, the first condition implies that 
there exists an open neighborhood $\U(\nu_Y)$ of $\ov\M_{g,\f}(Y,\be_Y;J)$
in $\X_{g,\f}(X,\be)$ such that 
$$\ov\M_{g,\f}(X,\be;J,\nu_Y)\cap \U(\nu_Y) = 
\ov\M_{g,\f}(Y,\be_Y;J,\nu_Y).$$
In addition, choose a multi-section $\vph$ of the bundle 
$\bL_{\f,\a}\lra\X_{g,\f}(X,\be)$ so that $\vph$ is transverse to the zero set
on every stratum of $\ov\M_{g,\f}(Y,\be_Y;J,\nu_Y)$ and every stratum of
$\ov\M_{g,\f}(X,\be;J,\nu)$.
This implies that the dimension of every stratum of $\ov\M_{g,\f}^{\vph}(Y,\be_Y;J,\nu_Y)$
is at most the rank~\e_ref{cokrk_e} of the bundle $\cok(D_J^{NY})$ over 
$\ov\M_{g,\f}(Y,\be_Y;J)$ and the equality holds only for the main stratum.\\

\noindent
By the middle assumption on $\nu_Y$ above,
for every element $[u]$ of $\ov\M_{g,S}(Y,\be_Y;J,\nu_Y)$
the linearization
$$D_{J,\nu_Y;u}^X\!: 
\H_u\oplus L^p_1(\Si_u;u^*TX)\lra L^p(\Si_u;T^*\Si_u^{0,1}\!\otimes\!_{\C}u^*TX)$$
of the section $\dbar_J\!+\!\nu_Y$ for maps to $X$ restricts to
the linearization
$$D_{J,\nu_Y;u}^Y\!: 
\H_u\oplus L^p_1(\Si_u;u^*TY)\lra L^p(\Si_u;T^*\Si_u^{0,1}\!\otimes\!_{\C}u^*TY)$$
of the section $\dbar_J\!+\!\nu_Y$ for maps to $Y$.
Thus, $D_{J,\nu_Y;u}^X$ descends to a Fredholm operator
$$D_{J,\nu_Y;u}^{NY}\!: L^p_1(\Si_u;u^*NY)\lra L^p(\Si_u;T^*\Si_u^{0,1}\!\otimes\!_{\C}u^*NY).$$
If $\nu_Y$ is sufficiently small,  by the last assumption in Theorem~\ref{main_thm} 
the operator 
\begin{equation*}\begin{split}
D_{J,\nu_Y,\vr;\u}^{NY}\!\equiv\!\big(D_{J,\nu_Y;u}^{NY}\big)_{\vr}\!: 
\big\{\xi\!\in\!L^p_1(\Si_u;u^*NY)\!:\,\xi(z_j(u))\!\in\!\Im\,d_{w_j}^{NY}f_j
~\forall\,j\!\in\!S\big\} \qquad\qquad&\\
\lra L^p(\Si_u;T^*\Si_u^{0,1}\!\otimes\!_{\C}u^*NY)&
\end{split}\end{equation*}
is injective for every $[\u]\!\in\!\ov\M_{g,\f}(Y,\be_Y;J,\nu_Y)$
as in~\e_ref{udfn_e}.
Thus, the cokernels of these operators still form an oriented vector orbi-bundle
over $\ov\M_{g,\f}(Y,\be_Y;J,\nu_Y)$ of rank~\e_ref{cokrk_e}, which will be denoted 
by $\cok(D_{J,\nu_Y,\vr}^{NY})$.
Furthermore, $\ov\M_{g,\f}(Y,\be_Y;J,\nu_Y)$ is compact
(because $\ov\M_{g,\f}(Y,\be_Y;J)$ is) and is a union of connected components of 
 $\ov\M_{g,\f}(X,\be;J,\nu_Y)$ by Proposition~\ref{horreg_prp}.\\

\noindent
The left-hand side of \e_ref{mainthm_e2} is the number of elements of 
$$\ov\M_{g,\f}^{\vph}(X,\be;J,\nu_Y\!+\!\nu_X)
\subset \X_{g,S}(X,\be)\times\prod_{j\in S}M_j$$
that lie in a small neighborhood of $\ov\M_{g,\f}^{\vph}(Y,\be_Y;J,\nu_Y)$
for any sufficiently small and generic~$\nu_X$.
The map component of any such element must be of the form $\exp_{u_{\ups}}\!\xi$, where
\begin{itemize}
\item $([u],(w_j)_{j\in S})\!\in\!\ov\M_{g,\f}(Y,\be_Y;J,\nu_Y)$
is an element of a fixed stratum, i.e.~the topological structure of~$\Si_u$ is fixed;
\item $\ups$ is a small gluing parameter for $\Si_u$ consisting of the smoothings
of the nodes of~$\Si_u$;
\item $u_{\ups}\!:\Si_{u_{\ups}}\!\lra\!Y$ is the approximately $(J,\nu_Y)$-map
corresponding to~$\ups$  as in \cite[Section~3]{gluing};
\item $\xi\!\in\!L^p_1(\Si_{\ups};u_{\ups}^*TX)$ is small with respect to
the $\|\cdot\|_{\ups,p,1}$-norm of \cite[Section~3]{LiT} and satisfies
\BE{XvsY_e0}\begin{split}
& \{\dbar_J\!+\!\nu_Y\}u_{\ups} 
+ D_{J,\nu_Y;u_{\ups}}\xi+\nu_X(u_{\ups})+N_{\ups}(\xi)=0,\\
&\xi(z_j(u_{\ups}))\in\Im\,(d_{w_j}^{NY}f_j)+T_{f_j(w_j)}Y \quad\forall\,j\!\in\!S,
\end{split}\EE
where  $N_{\ups}$ is a combination of a term quadratic in $\xi$ and 
a term which is linear in $\xi$ and~$\nu_X$.\\
\end{itemize}

\noindent
Projecting \e_ref{XvsY_e0} to $NY$, we obtain
\BE{XvsY_e1}\begin{split}
& D_{J,\nu_Y;u_{\ups}}^{NY}\ze+\nu_X^{\perp}(u_{\ups})+N_{\ups}^{\perp}(\ze)=0,\\
& \ze\in L^p_1(\Si_{u_{\ups}};u_{\ups}^*NY),\quad
\ze(z_j(u_{\ups}))\in\Im\,(d_{w_j}^{NY}f_j) ~~~\forall\,j\!\in\!S. 
\end{split}\EE
This equation has no small solutions in $\vph^{-1}(0)$
away from the subset of elements 
$$\u\!\equiv\!([u],(w_j)_{j\in S})\in\ov\M_{g,\f}^{\vph}(Y,\be_Y;J,\nu_Y)$$
for which $\nu_X^{\perp}(\u)$ lies in the image of $D_{J,\nu_Y,\vr;\u}^{NY}$,
i.e.~the projection $\bar\nu_X(\u)$ of 
$\nu_X(\u)$ to $\cok(D_{J,\nu_Y,\vr;\u}^{NY})$ is zero.
For dimensional reasons, all zeros of $\bar\nu_X$ lie in the main stratum 
$$\M_{g,\f}^{\vph}(Y,\be_Y;J,\nu_Y)\equiv
\ov\M_{g,\f}^{\vph}(Y,\be_Y;J,\nu_Y)\cap \M_{g,\f}(Y,\be_Y;J,\nu_Y).$$
Thus, only $\M_{g,\f}^{\vph}(Y,\be_Y;J,\nu_Y)$ contributes
to the left-hand side in~\e_ref{mainthm_e2}.
In this case equation~\e_ref{XvsY_e1} no longer involves~$\ups$ and
thus $u_{\ups}\!=\!u$.
Since $\vph$ vanishes transversally on $\M_{g,\f}(Y,\be_Y;J,\nu_Y)$
and  $\bar\nu_X$ on $\M_{g,\f}^{\vph}(Y,\be_Y;J,\nu_Y)$,
the contribution of the main stratum to
the left-hand side is the signed cardinality of the oriented zero-dimensional orbifold
$$\M_{g,\f}^{\vph}(Y,\be_Y;J,\nu_Y)\cap\bar\nu_X^{-1}(0).$$
As~$\bar\nu_X$ extends to a continuous multi-section of the orbi-bundle
\BE{cokeuler_e0}\cok(D_{J,\nu_Y,\vr}^{NY})\lra \ov\M_{g,\f}^{\vph}(Y,\be_Y;J,\nu_Y),\EE
which is transverse to the zero set over every stratum, 
the left-hand side of~\e_ref{mainthm_e2} is the euler class of 
the bundle~\e_ref{cokeuler_e0}
evaluated on $\ov\M_{g,\f}^{\vph}(Y,\be_Y;J,\nu_Y)$.
While the operators $D_{J,\nu_Y;\u}^{NY}$ and $D_{J;u}^{NY}$ are not the same,
they are homotopic through operators keeping the dimension of the cokernels fixed
and thus define orbi-bundles with the same euler class, as needed.\\

\noindent
The above argument requires some notion of smoothness for the strata of
$\X_{\T,\f}(X)$ or at least $\X_{\T,\f}(Y)$.
If the domain curve $\Si_u$ of $[u]$ with its marked points is stable
for every element $([u],(w_j)_{j\in S})$ of $\ov\M_{g,\f}(Y,\be_Y;J)$,
then every stratum $\X_{\T,\f}(X)$ meeting $\ov\M_{g,\f}(Y,\be_Y;J)$
is a smooth Banach orbifold.
The topological aspects of the resulting setting are sorted out in~\cite{Mc0},
and the above argument suffices in such cases.
These include the cases of 
Theorem~\ref{FanoGV_thm} (with $2g\!+\!|S|\!\ge\!3$, which can be assumed) 
and Corollary~\ref{LeP_crl} (since the genus of $Y_{\al}$ is positive), 
but not of Example~\ref{CY_eg} or the specific cases 
of Examples~\ref{inY_eg} or~\ref{transversetoY_eg}.\\

\noindent
In general, $\X_{\T}(X)$ is a subspace of a product of main strata 
$\X_{g_i,S_i}^0(X,\be_i)$ for some $g_i$, $S_i$, and $\be_i$ and
the restriction of $\Ga^{0,1}_{g,S}(X,\be;J)$ is the direct sum 
of the pull-backs of the corresponding bundles over the components 
of the product.
If for every $([u],(w_j)_{j\in S})\!\in\!\ov\M_{g,\f}(Y,\be_Y;J)$
and every unstable component $\Si_{u;i}$ of $\Si_u$
the restriction of $u$ to $\Si_{u;i}$ is regular in the appropriate sense,
then $\nu$ can be taken to be a smooth section of the components of 
$\Ga^{0,1}_{g,S}(X,\be;J)$ coming from the ``stable parts" of~$\T$;
as in the previous paragraph there is a well-defined notion of smoothness
over these components.
This is done explicitly in \cite[Section~2]{RT2}.
The resulting extension of the previous paragraph then covers 
the specific cases of Examples~\ref{inY_eg} and~\ref{transversetoY_eg}.\\

\noindent
Finally, for an arbitrary symplectic manifold $(X,\om)$,
the ``notion" of smoothness is described by introducing smooth 
finite-dimensional approximations to $\ov\M_{g,S}(X,\be;J)$.
This is done in the next section.

\subsection{Virtual setting}
\label{mainpf_subs}

\noindent
Continuing with the notation of Section~\ref{config_subs},
we now recall the virtual fundamental class setup of~\cite{FuO} 
and~\cite{LiT} and then reformulate the argument of 
Section~\ref{semipos_subs} for \e_ref{mainthm_e2} in the general case.\\

\noindent
An \sf{atlas for $\ov\M_{g,S}(X,\be;J)$} is a collection
$\{(\U_{\al},E_{\al})\}_{\al\in\A}$, where 
\begin{itemize}
\item $\{\U_{\al}\}_{\al\in\A}$ is an open cover of
$\ov\M_{g,S}(X,\be;J)$ in $\X_{g,S}(X,\be)$  and 
$E_{\al}\!\subset\!\Ga_{g,S}^{0,1}(X,\be;J)|_{\U_{\al}}$
is a topological (finite-rank) vector orbi-bundle over~$\U_{\al}$;
\item $\dbar_J^{-1}(E_{\al})$ is a smooth orbifold and
$\dbar_J^{-1}(E_{\al})\cap\X_{T}(X)$ is a smooth sub-orbifold of $\dbar_J^{-1}(E_{\al})$
of the codimension corresponding to $\T$ (twice the number of nodes)
for every stratum $\X_{T}(X)$;
\item the restriction of $E_{\al}$ to $\dbar_J^{-1}(E_{\al})$
is a smooth vector orbi-bundle and the restriction of $\dbar_J$ to 
$\dbar_J^{-1}(E_{\al})$ 
is a smooth section of $E_{\al}|_{\dbar_J^{-1}(E_{\al})}$;
\item for every $[u]\!\in\!\ov\M_{g,S}(X,\be;J)\!\cap\!\dbar_J^{-1}(E_{\al})\!\cap\!
\dbar_J^{-1}(E_{\al'})$, there exists $\ga\!\in\!\A$ such that
$$[u]\in\U_{\ga}\subset\U_{\al}\cap\U_{\al'}\,,
\qquad E_{\al},E_{\al'}\big|_{U_{\ga}}\subset E_{\ga}\,,$$
the restrictions of $E_{\al}$ and $E_{\al'}$ to
$\dbar_J^{-1}(E_{\ga})\cap\X_{T}(X)$
are smooth orbifold subbundles of the restriction of~$E_{\ga}$, and
the restriction of $\dbar_J$ to $\dbar_J^{-1}(E_{\ga})\cap\X_{T}(X)$
is transverse to $E_{\al}$ and~$E_{\al'}$;
\item for every $[u]\!\in\!\ov\M_{g,S}(X,\be;J)$, 
\BE{chartcond_e}
\Ga^{0,1}(X,u;J)=
\big\{D_{J;u}\xi\!:\,\xi\!\in\!L^p_1(\Si_u;u^*TX)\big\}+\ti{E}_{\al}|_u\,,
\EE
where $\ti{E}_{\al}|_u\!\subset\!\ti\Ga^{0,1}_{\T}(X;J)|_u$ 
is the preimage of $E_{\al}|_u$ 
under the quotient map 
$$\ti\Ga^{0,1}_{\T}(X;J)|_u\!\lra\!\Ga^{0,1}_{g,S}(X,\be;J)|_{[u]}\,.$$
\end{itemize}
Such collections $\{(\U_{\al},E_{\al})\}_{\al\in\A}$ 
are described in \cite[Section~12]{FuO} and \cite[Section~3]{LiT}.
An atlas for $\ov\M_{g,\f}(X,\be;J)$ is defined similarly,
with the domain of $D_{J;u}$ in \e_ref{chartcond_e} replaced by
$$\big\{\xi\!\in\!L^p_1(\Si_u;u^*TX)\!:\,
\xi(z_j(u))\in\Im\,d_{w_j}f_j~\forall\,j\!\in\!S\big\}$$
for an element $([u],(w_j)_{j\in S})$ of $\ov\M_{g,\f}(X,\be;J)$.
Such an atlas induces a compatible atlas for the total space 
of the restriction of the bundle $\bL_{\a,\f}$ to $\ov\M_{g,\f}(X,\be;J)$.\\

\noindent
A \sf{multi-section $\nu$ of $\Ga^{0,1}_{g,\f}(X,\be;J)$ for an atlas $\{(\U_{\al},E_{\al})\}_{\al\in\A}$} is a continuous multi-section
such that the restriction of $\nu$ to $\dbar_J^{-1}(E_{\al})$
is a smooth section of~$E_{\al}$.
Similarly, a {multi-section $\vph$ of $\bL_{\a,\f}$ for 
$\{(\U_{\al},E_{\al})\}_{\al\in\A}$} is a continuous multi-section
such that the restriction of $\vph$ to $\dbar_J^{-1}(E_{\al})$ is smooth.
A multi-section $\nu$ as above is \sf{regular} if the restriction of $\nu$
to $\dbar_J^{-1}(E_{\al})\!\cap\!\X_{\T,\f}(X)$ is transverse to 
the zero set in $E_{\al}$ for every $\al$ and~$\T$.
If $(\{(\U_{\al},E_{\al})\}_{\al\in A},\nu)$ is regular, 
$\ov\M_{g,\f}(X,\be;J,\nu)$ is stratified by smooth 
branched orbifolds of even dimensions.
The existence of regular multi-sections for a refinement of a subatlas 
is the subject of \cite[Chapter~1]{FuO} and \cite[Section~4]{Mc1}.\footnote{It is also
shown in \cite{FuO} and \cite{Mc1} that a regular multi-section $\nu$ 
determines a rational homology class;
however, this notion of virtual fundamental class is not necessary for 
defining GW-invariants or comparing the two sides of~\e_ref{mainthm_e2}.}
If $\nu$ is sufficiently small and regular and $\vph$ is generic,
the left-hand side of~\e_ref{mainthm_e2} is again the weighted number of 
elements of
$$\ov\M_{g,\f}^{\vph}(X,\be;J,\nu)\equiv \ov\M_{g,\f}(X,\be;J,\nu)\cap\vph^{-1}(0)$$
that lie in a small neighborhood of 
$$\ov\M_{g,\f}^{\vph}(Y,\be_Y;J)\equiv  \ov\M_{g,\f}(Y,\be_Y;J)\cap\vph^{-1}(0)$$
in $\X_{g,\f}(X,\be)$.\\

\noindent
By \cite[Chapter~3]{FuO} and \cite[Section~3]{LiT}, 
pairs $(\U_{Y;\al},E_{Y;\al})$ for an atlas for 
$$\ov\M_{g,S}(Y,\be_Y;J)\times\prod_{j\in S}\!M_j$$ 
that restrict to an atlas for $\ov\M_{g,\f}(Y,\be_Y;J)$
can be obtained in the following way.
Given $\u\!=\!([u],(w_j)_{j\in S})$, choose
\begin{itemize}
\item a neighborhood $V_{Y;u}$ of $u(\Si_u)$ in $Y$;
\item a representative $u\!:\Si_u\!\lra\!Y$ for $[u]$; 
\item universal family of deformations $\W_u\!\lra\!\De_u$ of $\Si_u$
with its marked points (thus $\Si_u\!\subset\!\W_u$);
\item  a finite-dimensional subspace
$$\cE_{Y;\u}\subset \Ga_c\big(\W_u^*\!\times\!V_{Y;u};
\pi_1^*(T^*\W_u^v)^{0,1}\!\otimes_{\C}\!\pi_2^*TY\big),$$
where $\W_u^*\!\subset\!\W_u$ is the subspace of smooth points of the fibers,
$T\W_u^v\!\subset\!T\W_u$ is the vertical tangent space,
and $\Ga_c$ denotes the space of
smooth compactly supported bundle sections, such that
$$\Ga(\Si_u;T^*\Si_u^{0,1}\!\otimes_{\C}\!u^*TY)=
\big\{D_u\xi\!: \xi\!\in\!\Ga(\Si_u;u^*TY),\,\xi(z_i(u))\!\in\!\Im\,d_{w_j}f_j\,
\forall\,j\!\in\!S\big\}
+\{\id\!\times\!u\}^*\cE_{Y;\u}\,$$
if $\u\!\in\!\ov\M_{g,\f}(Y,\be_Y;J)$; 
if $\u\!\not\in\!\ov\M_{g,\f}(Y,\be_Y;J)$, the point-wise condition on $\xi$ is omitted. 
\end{itemize}
If $\u'\!=\!([u'],(w_j')_{j\in S})$ with
$[u']\!\in\!\X_{g,S}(V_{Y;u},\be_Y)$ and $\Si_{u'}\!\in\!\De_u$, let 
$$\ti{E}_{Y;\u}|_{\u'}=\{\id\!\times\!u'\}^*\cE_{Y;\u}\,.$$
By \cite[Chapter~3]{FuO} and \cite[Section~3]{LiT},
$\U_{Y;\al}$ can be taken to be the image of a sufficiently small
neighborhood $\ti\U_{Y;\al}$ of $\u$ in the space of $L^p_1$-maps
from the fibers of $\W_u\!\lra\!\De_u$ to $X$ under the equivalence 
relation and $E_{Y;\al}$ the image of the bundle formed by the spaces 
$\ti{E}_{Y;\u}|_{\u'}$ over $\ti\U_{Y;\al}$.
With these choices, $\dbar_J^{-1}(E_{Y;\al})$ consists of equivalence classes
of smooth maps to~$Y$.\\

\noindent
Fix a $J$-regularized tubular neighborhood $(\pi_Y\!:U_Y\!\lra\!Y,TU_Y^h)$ of $Y$
in~$X$.
Using $\pi_Y$ and $d\pi_Y|_{TU_Y^h}^{-1}$, each $\cE_{Y;\u}$ can be extended to 
a finite-dimensional subspace 
$$\cE_{X|Y;\u}\subset \Ga_c\big(\W_u^*\!\times\!V_{X;u};
\pi_1^*(T^*\W_u^v)^{0,1}\!\otimes_{\C}\!\pi_2^*TU_Y^h\big)
\subset \Ga_c\big(\W_u^*\!\times\!V_{X;u};
\pi_1^*(T^*\W_u^v)^{0,1}\!\otimes_{\C}\!\pi_2^*TX\big)$$
for a neighborhood $V_{X;u}$ of $V_{Y;u}$ in $U_Y\!\subset\!Y$.
A larger subspace
$$\cE_{X;\u}
\subset \Ga_c\big(\W_u^*\!\times\!V_{X;u};
\pi_1^*(T^*\W_u^v)^{0,1}\!\otimes_{\C}\!\pi_2^*TX\big)$$
can then be chosen so that 
$$\Ga(\Si_u;T^*\Si_u^{0,1}\!\otimes_{\C}\!u^*TX)=
\big\{D_u\xi\!: \xi\!\in\!\Ga(\Si_u;u^*TX),\,\xi(z_i(u))\!\in\!\Im\,d_{w_j}f_j\,
\forall\,j\!\in\!S\big\}
+\{\id\!\times\!u\}^*\cE_{X;\u}\,,$$
whenever $[u]\!\in\!\ov\M_{g,\f}(Y,\be_Y;J)$. 
This gives rise to a pair $(\U_{X;\al},E_{X;\al})$ for an atlas for
$\ov\M_{g,\f}(X,\be;J)$; the union of such pairs covers $\ov\M_{g,\f}(Y,\be_Y;J)$.
Since $\ov\M_{g,\f}(Y,\be_Y;J)$ is a union of components of $\ov\M_{g,\f}(X,\be;J)$,
this sub-collection of an atlas is sufficient for determining the left-hand side
of~\e_ref{mainthm_e2}.
Similarly, using $\pi_Y$, $d\pi_Y|_{TU_Y^h}^{-1}$, and a bump function around $Y$ 
with support in~$U_Y$, any multi-section of 
$$\pi_1^*\Ga_{g,S}^{0,1}(Y,\be_Y;J)\lra \X_{g,S}(Y,\be_Y)
\times\prod_{j\in S}\!M_j$$
for the atlas $(\{(\U_{Y;\al},E_{Y;\al})\}_{\al\in\A})$ gives rise to 
a multi-section $\nu$ of 
$$\Ga_{g,\f}^{0,1}(X,\be;J)\lra \X_{g,\f}(X,\be)$$
for the atlas $(\{(\U_{X;\al},E_{X;\al})\}_{\al\in\A})$
such that for every element $[\u]\!\in\!\X_{g,\f}(X,\be)$
$$\nu([\u])\in 
L^p\big(\Si_u;T^*\Si_u^{0,1}\!\otimes\!u^*TU_Y^h\big)$$
for every $[\u]\!=\!([u],(w_j)_{j\in S})\!\in\!\X_{g,\f}(U_Y,\be_Y)$.\\

\noindent
Let $\nu\!=\!\nu_Y\!+\!\nu_X$ be a regular multi-section of $\Ga^{0,1}_{g,\f}(X,\be)$
for atlas for $\ov\M_{g,\f}(X,\be;J)$ as above so that 
\begin{itemize}
\item for every $\u\!=\!([u],(w_j)_{j\in S})\!\in\!\X_{g,\f}(X,\be)$
with $[u]\!\in\!\X_{g,S}(U_Y,\be_Y)$,
$$\nu_Y(\u)\in L^p(\Si_u;T^*\Si_u^{0,1}\!\otimes_{\C}\!u^*TU_Y^h);$$
\item $\nu_Y|_{\X_{g,\f}(Y,\be_Y)}$ is a regular multi-section 
of $\Ga^{0,1}_{g,\f}(Y,\be_Y)$
so that $\ov\M_{g,\f}(Y,\be_Y;J,\nu_Y)$ is stratified by smooth 
branched orbifolds of the expected dimensions and the dimension of the main stratum
$$\M_{g,\f}(Y,\be_Y;J,\nu_Y) \equiv \ov\M_{g,\f}(Y,\be_Y;J,\nu_Y)
\cap \bigg(\X_{g,S}^0(Y,\be_Y)\times\prod_{j\in S}\!M_j\bigg)$$
is $\dim_{g,\f}(Y,\be_Y)$;
\item $\nu_X$ is small relative to $\nu_Y$. 
\end{itemize} 
The previous paragraph implies that such multi-sections $\nu_Y$ exist.
By Proposition~\ref{horreg_prp}, the first condition implies that 
there exists an open neighborhood $\U(\nu_Y)$ of $\ov\M_{g,\f}(Y,\be_Y;J)$
in $\X_{g,\f}(X,\be)$ such that 
$$\ov\M_{g,\f}(X,\be;J,\nu_Y)\cap \U(\nu_Y) = 
\ov\M_{g,\f}(Y,\be_Y;J,\nu_Y).$$
In addition, choose a multi-section $\vph$ of the bundle 
$\bL_{\f,\a}\lra\X_{g,\f}(X,\be)$ for the above atlas 
so that $\vph$ is transverse to the zero set
on every stratum of $\ov\M_{g,\f}(Y,\be_Y;J,\nu_Y)$ and every stratum of
$\ov\M_{g,\f}(X,\be;J,\nu)$.\\

\noindent
For each $\al\!\in\!\A$ and 
$\u\!\in\!\ov\M_{g,\f}(Y,\be_Y;J,\nu_Y)\!\cap\!\U_{Y;\al}$, let
$$\cD_{\nu_Y,\al;\u}\!: T_{\u}\dbar_J^{-1}(E_{X;\al}) \lra E_{X;\al}$$
be the linearization of the section $\dbar_J\!+\!\nu_Y$ over
$\dbar_J^{-1}(E_{X;\al})$ along the zero set.
The kernel of $\cD_{\nu_Y,\al;\u}$ is the tangent space of 
$\ov\M_{g,\f}(Y,\be_Y;J,\nu_Y)$ at~$\u$.
If $\al$ and $\ga$ are as in the overlap condition in the definition of
an atlas above, then
\begin{gather*}
E_{X;\al}\cap \Im\,\cD_{\nu_Y,\ga;\u}=\Im\,\cD_{\nu_Y,\al;\u}
\quad\forall\,\u\in\ov\M_{g,\f}(Y,\be_Y;J,\nu_Y)\cap\U_{Y;\ga}\,,\\
\dim\dbar_J^{-1}(E_{X;\ga})-\dim\dbar_J^{-1}(E_{X;\al})
=\rk\, E_{X;\ga}-\rk\, E_{X;\al}.
\end{gather*}
Thus, the inclusion $T\dbar_J^{-1}(E_{X;\al})\lra T\dbar_J^{-1}(E_{X;\ga})$
over $\ov\M_{g,\f}(Y,\be_Y;J,\nu_Y)\!\cap\!\U_{Y;\ga}$ induces isomorphisms
$$\cok(\cD_{\nu_Y,\al;\u})\lra \cok(\cD_{\nu_Y,\ga;\u}).$$
It follows that these vector spaces form an orbi-bundle $\cok(\cD_{\nu_Y})$
over $\ov\M_{g,\f}(Y,\be_Y;J,\nu_Y)$.
By the last requirement in the definition of an atlas
and condition~(b) in Theorem~\ref{main_thm}, the homomorphism
$$\cok(\cD_{\nu_Y,\al;\u})\lra\cok(D_{J;\u}^{NY})$$
induced by the inclusion $E_{X;\al}\lra\Ga^{0,1}_{g,\f}(X,\be;J)$ 
followed by the projections to $NY$ and the cokernel is 
surjective for all $\u\!\in\!\ov\M_{g,\f}(Y,\be_Y;J,\nu_Y)\!\cap\!\U_{Y;\al}$,
if $\nu_Y$ is sufficiently small.
A dimension count then shows that this homomorphism is an isomorphism 
(the injectivity also follows from Proposition~\ref{horreg_prp}).
Thus, the orbi-bundles
$$\cok(\cD_{\nu_Y}),\cok(D_J^{NY})\lra\ov\M_{g,\f}(Y,\be_Y;J,\nu_Y) $$
are isomorphic.\\

\noindent
The left-hand side of \e_ref{mainthm_e2} is the number of elements of 
$$\ov\M_{g,\f}^{\vph}(X,\be;J,\nu_Y\!+\!\nu_X)
\subset \X_{g,S}(X,\be)\times\prod_{j\in S}\!M_j$$
that lie in a small neighborhood of $\ov\M_{g,\f}^{\vph}(Y,\be_Y;J,\nu_Y)$
for a small generic multi-section $\nu_X$.
The number of such elements near 
$\ov\M_{g,\f}^{\vph}(Y,\be_Y;J,\nu_Y)\!\cap\!\U_{Y;\al}$
is the number of solutions of
$$  \cD_{\nu_Y,\al;\u}\xi+\nu_X(\u)+N_{\al}(\xi)=0,
\qquad \xi\in T_{\u}\dbar_J^{-1}(E_{X;\al}), $$
with small $\xi$, where
 $N_{\al}$ is a combination of a term quadratic in $\xi$ and 
a term which is linear in $\xi$ and~$\nu_X$.
This equation has no solutions in $\vph^{-1}(0)$
away from the subset of elements 
$$\u\in\ov\M_{g,\f}^{\vph}(Y,\be_Y;J,\nu_Y)$$
for which $\nu_X(\u)$ lies in the image of $\cD_{\nu_Y,\al;\u}$,
i.e.~the projection $\bar\nu_X(\u)$ to $\cok(\cD_{\nu_Y,\al;\u})$ is zero.
Since $\vph$ vanishes transversally on $\ov\M_{g,\f}(Y,\be_Y;J,\nu_Y)$
and $\bar\nu_X$ on $\ov\M_{g,\f}^{\vph}(Y,\be_Y;J,\nu_Y)$,
the left-hand side of \e_ref{mainthm_e2} is the signed cardinality of 
oriented zero-dimensional orbifold
$$\ov\M_{g,\f}^{\vph}(Y,\be_Y;J,\nu_Y)\cap\vph^{-1}(0).$$
By the definition, this is also the euler class of $\cok(\cD_{\nu_Y})$
evaluated on $\ov\M_{g,\f}^{\vph}(Y,\be_Y;J,\nu_Y)$,
which by the above isomorphism of cokernel bundles equals to the right-hand side of \e_ref{mainthm_e2}.\\

\vspace{.2in}

\noindent
{\it Department of Mathematics, SUNY Stony Brook, Stony Brook, NY 11790-3651\\
azinger@math.sunysb.edu}\\


\begin{thebibliography}{99}


\bibitem{ACGH} E.~Arbarello, M.~Cornalba, P~Griffiths, J.~Harris, 
Geometry of Algebraic Curves, Vol.~I, 
Springer-Verlag, New York, 1985.

\bibitem{AsMo} P.~Aspinwall and D.~Morrison, 
{\it Topological field theory and rational curves}, 
Comm. Math.~Phys. 151 (1993), 245--262.

\bibitem{BP} J.~Bryan and R.~Pandharipande, 
{\it BPS states of curves in Calabi-Yau $3$-folds}, 
Geom.~Top.~5 (2001), 287-318.

\bibitem{FHoS} A.~Floer, H.~Hofer, and D.~Salamon, 
{\it Transversality in elliptic Morse theory for the symplectic action},
Duke Math.~J.~80 (1996), no.~1, 251-292.

\bibitem{FuO} K.~Fukaya and K.~Ono,
{\it Arnold Conjecture and Gromov-Witten invariant},
Topology 38 (1999), no.~5, 933--1048.

\bibitem{GV1} R.~Gopakumar and C.~Vafa, 
{\it M-theory and topological Strings~I}, hep-th/9809187.

\bibitem{GV2} R.~Gopakumar and C.~Vafa,
{\it M-theory and topological strings~II}, hep-th/9812127.

\bibitem{HuLtR} J.~Hu, T.-J.~Li, and Y.~Ruan, 
{\it Birational cobordism invariance of uniruled symplectic manifolds},
Invent.~Math.~172  (2008),  no.~2, 231--275.

\bibitem{HuR} J.~Hu and Y.~Ruan, 
{\it Positive divisors in symplectic geometry}, math/0802.0590.

\bibitem{IvSh} S.~Ivashkovich and V.~Shevchishin,
{\it Pseudo-holomorphic curves and envelopes of meromorphy of two-spheres in 
$\C P^2$}, math.CV/9804014.

\bibitem{KKO} B.~Kim, A.~Kresch, and Y.-G.~Oh, 
{\it A compactification of the space of maps from curves}, preprint.

\bibitem{KLi} Y.-H.~Kiem and J.~Li,
{\it Gromov-Witten invariants of varieties with holomorphic 2-forms},
math/0707.2986.

\bibitem{KlP} A.~Klemm and R.~Pandharipande,
{\it Enumerative geometry of Calabi-Yau 4-folds},
Comm.~Math.~Phys.~261 (2006), no.~2, 451--516.

\bibitem{KoM} M.~Kontsevich and Yu.~Manin,
{\it Gromov-Witten classes, quantum cohomology, and enumerative geometry}, 
Comm.~Math.~Phys.~164 (1994), no.~3, 525--562. 

\bibitem{LiT}  J.~Li and G.~Tian, 
{\it Virtual moduli cycles and Gromov-Witten invariants of general symplectic manifolds}, 
Topics in Symplectic \hbox{$4$-Manifolds},
Internat.~Press, 1998.

\bibitem{LiZ}  J.~Li and A.~Zinger, 
{\it On Gromov-Witten invariants of a quintic threefold
and a rigidity conjecture}, Pacific J.~Math.~233 (2007), no.~2, 417-480.

\bibitem{LeP} J.~Lee and T.~Parker, 
{\it A structure theorem for the Gromov-Witten invariants of Kahler surfaces}, 
JDG  77  (2007),  no.~3, 483--513.

\bibitem{LtR} T.-J.~Li and Y.~Ruan, 
{\it Uniruled symplectic divisors}, math/0711.4254. 

\bibitem{MaP}  D.~Maulik and R.~Pandharipande, 
{\it New calculations in Gromov--Witten theory},
math.AG/0601395.

\bibitem{Mc0} D.~McDuff, {\it The virtual moduli cycle}, 
Northern California Symplectic Geometry Seminar, 
pp 73--102, AMS Transl.~Ser.~2, 196, 1999.

\bibitem{Mc1} D.~McDuff,
{\it Groupoids, branched manifolds, and multisections},
J.~Symplectic Geom.~4 (2006), no.~259--315.


\bibitem{Mc} D.~McDuff,
{\it Hamiltonian $S^1$-manifolds are uniruled}, 
Duke Math.~J.~146 (2009), no.~3, 449--507.



\bibitem{McS} D.~McDuff and D.~Salamon,
{\it $J$-holomorphic Curves and Symplectic Topology},
AMS~2004.

\bibitem{McTo} D.~McDuff and S.~Tolman,
{\it Topological properties of Hamiltonian circle actions}, 
IMRP Int.~Math.~Res.~Pap (2006), 1--77.

\bibitem{MirSym} K.~Hori, S.~Katz, A.~Klemm, R.~Pandharipande,
R.~Thomas, C.~Vafa, R.~Vakil, and E.~Zaslow, {\it Mirror Symmetry},
Clay Math.\ Inst., Amer.\ Math.\ Soc.,~2003. 

\bibitem{P1} R.~Pandharipande,
{\it Hodge integrals and degenerate contributions}, 
Comm.~Math.~Phys.~208  (1999),  no.~2, 489--506.

\bibitem{P2} R.~Pandharipande,
{\it Three questions in Gromov-Witten theory},
Proceedings of ICM, Beijing (2002),  503--512.
 
\bibitem{RT}  Y.~Ruan and G.~Tian, {\it A mathematical theory of quantum cohomology},  
JDG 42 (1995),  no.~2, 259--367.

\bibitem{RT2}  Y.~Ruan and G.~Tian, {\it Higher genus symplectic invariants 
and sigma models coupled with gravity},  
Invent.~Math.~130  (1997),  no.~3, 455--516. 

\bibitem{SeSi} R.~Seeley and I.~Singer,
{\it Extending $\dbar$ to singular Riemann surfaces},
J.~Geom.~Phys.~4 (1988), no~1, 121--136. 

\bibitem{Si} B.~Siebert, {\it Gromov-Witten invariants for general
symplectic manifolds}, dga-ga/9608005.

\bibitem{Sh} V.~Shevchishin, 
{\it Pseudoholomorphic curves and the symplectic isotopy problem},
math/0010262.

\bibitem{Taubes} C.~Taubes, 
{\it Counting pseudo-holomorphic submanifolds in dimension 4},
Seiberg-Witten and Gromov invariants for symplectic 4-manifolds,  99--161, 
First Int.~Press Lect.~Ser~2. 


\bibitem{Voisin} C.~Voisin, {\it A mathematical proof of a formula of
Aspinwall and Morrison}, Comp.~Math.~104 (1996), no.~2, 135--151.

\zr{gluing} 
{\it Enumerative vs.~symplectic invariants and obstruction bundles},
J.~Symplectic Geom.~2 (2004), no.~4, 445--543.

\zr{g1comp2} 
{\it Reduced genus-one Gromov-Witten invariants}, math.SG/0507103.

\zr{divisorGWs}
{\it On relative Gromov-Witten invariants in genus zero}, 
in preparation.

\end{thebibliography}
\end{document}